%% file: main.tex
\documentclass[hidelinks,onefignum,onetabnum]{siamart220329}
\usepackage{ulem}

\input{ex_shared}

\ifpdf
\hypersetup{
  pdftitle={WAN discretization of PDEs: best approximation, stabilization and essential boundary conditions},
  pdfauthor={Silvia Bertoluzza, Erik Burman, Cuiyu He}
}
\fi

\newcommand{\aform}{\mathcal{A}}
\newcommand{\ff}{\mathcal{F}}
\newcommand{\parameterT}{\mathcal{P}_{\bs\theta}}
\newcommand{\parameterZ}{\mathcal{P}_{\bs\eta}}
\newcommand{\bs}[1]{\boldsymbol{#1}}
\newcommand{\azang}{a}
\newcommand{\bzang}{b}
\newcommand{\Azang}{A}
\newcommand{\Bzang}{B}
\newcommand{\tnorm}[1]{\vvvert #1 \vvvert}
\newcommand{\wstar}{u^*_{\bs\theta}}
\newcommand{\minimum}{\sigma^*}
\newcommand{\tw}{\widetilde w_{\bs \theta}}

\newtheorem{lem}{Lemma}
\newtheorem{ex}{Example}
\newtheorem{assumption}[theorem]{Assumption}

\begin{document}

\maketitle

\begin{abstract} 
In this paper, we provide a theoretical analysis of the recently introduced weakly adversarial networks (WAN) method, used to approximate partial differential equations in high dimensions. We address the existence and stability of the solution, as well as approximation bounds. We also propose two new stabilized WAN-based formulas that avoid the need for direct normalization. Furthermore, we analyze the method's effectiveness for the Dirichlet boundary problem that employs the implicit representation of the geometry. We also devise a pseudo-time XNODE neural network for static PDE problems, yielding significantly faster convergence results than the classical DNN.
\end{abstract}



\section{Introduction}
Recently there has been a vast interest in approximating partial differential equations (PDE) using neural networks and machine learning techniques. In this note, we will consider the weak adversarial networks (WAN) method introduced by Zang et al. \cite{ZBYZ20}. 
The idea is to rewrite the weak form of the PDE as a saddle point problem whose solution is obtained by approximating both the trial (primal) and the test (adversarial) space through neural networks. In \cite{ZBYZ20}, the method was tested on various PDEs, tackling different challenging issues such as high dimension, nonlinearity, and nonconvexity of the domain. It was subsequently applied for the inverse problems in high dimension \cite{BYZZ20} and for the parabolic problems \cite{oliva2021fast}, with quite promising results. 

However,  as often happens for neural network methods for numerical PDEs, rigorous theoretical results on the capability of WANs to approximate the solution of a given PDE still need to be improved. The most critical issues must be addressed are  
 the discrete solution's existence, stability, and approximation properties. Due to the inherent nature of neural network function classes, even the issue of the existence of a discrete solution is far from a trivial one.
Indeed, fixed architecture neural network classes are generally neither convex nor closed \cite{petersen2021topological,mahan2021nonclosedness}. Therefore, a global minimum for a cost functional in one of such classes might not exist. Unsurprisingly, as we are ultimately dealing with a saddle point problem, a suitable choice of the test (adversarial) network class will play a vital role in the analysis. The lack of linearity of the trial (primal) network class will imply the need for a strengthened inf-sup condition (see \cref{inf-sup} in the following), which, however, will not, in general, be enough to guarantee the existence and uniqueness of a global minimizer. Indeed, due to the non-closedness of neural network classes, it might not be possible to attain the minimum with an element belonging to the class. 

What we can prove, under suitable assumptions (see \cref{assumptionZT-1,inf-sup}),
 in a general abstract framework, is (1) the existence of at least one weakly converging minimizing sequence for the WAN cost functional, and,  (2) that all weak limits of weakly converging minimizing sequences satisfy a quasi-optimality bound similar to C\'ea's Lemma.
 More importantly, we further prove that a similar approximation bound will hold for the elements of the minimizing sequences sufficiently close to convergence. Combined with approximation bounds by the deep neural networks \cite{guhring2020expressivity}, this will guarantee that the WAN can, in principle, provide an arbitrarily good approximation to the continuous PDE solution.
Another crucial issue relates to the convergence of the optimization scheme used to solve the minimization problem.
Also this task is made difficult by the inherent topological properties of neural network classes.
 It is worth mentioning (see  \cite{petersen2021topological}) that the function class of Deep Neural Networks (DNN) lacks inverse stability in the $L^{p}$ and $W^{s,p}$ norms. In simple terms, the norm of the elements of the DNN function class does not control the norm of the associated parameter vector. As the optimization schemes indirectly act on the function class through the parameter space, this will negatively affect the minimization process. In particular, when, 
 the weak limit of the minimizing sequence does not belong to the function class, it can be proved that the sequence of the Euclidean norms of the corresponding parameter vectors explodes \cite{petersen2021topological}. 
 
In the  WAN framework, aforementioned problems are integrated with the problems related to the inexact evaluation of the cost functional, which is defined as a supremum over the elements of the adversarial network and requires solving an optimization problem that, for the classical WAN method, become ill-posed due to the presence of direct normalization, and is therefore subject to a possibly relevant error. If this error becomes comparable, or even dominant, compared to the value of the cost functional itself, the overall optimization procedure will lose effectiveness and likely display oscillations. To mitigate this phenomenon, developing more stable and accurate methods for evaluating the operator norm is crucial. In the framework of WAN, we propose two alternative ways of evaluating the operator, that avoid direct normalization and improve the overall convergence of the minimization procedure.

We then exploit the results 
for the second-order elliptic PDEs with essential boundary conditions. These are notoriously challenging 
as the construction of neural networks exactly vanishing on the boundary of a domain is extremely difficult, if not impossible. On the other hand, standard techniques, such as Nitsche's method, that impose Dirichlet boundary conditions weakly, rely on inverse inequalities that do not generally hold in the neural network framework. 
Adapting a strategy introduced, for finite elements, in \cite{DL20}, we propose to approximate the test space $H^1_0(\Omega)$ with a class of functions obtained by multiplying the elements of a given neural network class with a level set type weight, thus strongly enforcing the homogeneous boundary conditions on the test function class. Non-homogeneous boundary conditions are then imposed by penalization with a suitable boundary norm. We can show that the resulting discrete schemes fall in our abstract setting,  thus obtaining C\'ea's Lemma type quasi-optimal $H^1$ error bounds.

As the architecture of neural networks plays a crucial role in their performance, we test the newly proposed methods on different function classes of various structures. In particular, besides DNN, we focus 
 on residual-related networks, whose usage  \cite{he2016deep} was initially proposed to enhance image processing capabilities. These networks have also found application in various domains, including numerical PDEs \cite{karumuri2020simulator, zeng2022adaptive}. In a recent work by Oliva et al. \cite{oliva2021fast}, the XNODE network is proposed to solve parabolic equations. Numerical experiments have demonstrated that, compared to classical DNN networks, XNODE, 
 can substantially reduce the number of iterations required for optimization. 
 {This rapid convergence can be attributed to the structure of the XNODE model, which emulates a residual network, and to the direct incorporation of the initial condition into the model.} Besides testing our framework with XNODE architecture on a parabolic problem, we also introduce a new variant of the XNODE network, which we refer to as the pseudo-time XNODE method for stationary problems. Remarkably fast convergence is observed in the numerical results, even for nonlinear and high-dimensional static elliptic PDEs.
 
 The paper is organized as follows. 
In section \ref{sec:WAN_analysis}, we prove quasi-best approximation results and in Section \ref{sec:WAN_alternative} we
 propose two more stable equivalent formulations. In section \ref{sec:Dirichlet}, we  leverage our approach to allow for Dirichlet boundary conditions. Finally, the numerical results are provided in \cref{sec:num}.
We devote the remaining part of this  section to discuss the standard WAN in an abstract setting. Our framework covers a large class of problems without symmetry or coercivity assumptions, allowing for standard well-posed problems and certain non-standard data assimilation problems. We also cover a very general class of discretization spaces: while we have in mind neural networks, the only a priori assumptions that we make on our trial and test function classes is that they are function sets containing the identically vanishing function so that our results potentially applies to a much wider range of methods, provided that the inf-sup conditions \cref{assumptionZT-1,inf-sup} hold.

Throughout the paper, we assume that all forms, linear and bilinear, are evaluated exactly and that the resulting nonlinear optimization problems can be solved with sufficient accuracy. Needless to say, these problems are crucial for the actual performance of the method. Nevertheless, the quasi-best approximation results proved herein are a cornerstone for its reliability.

\subsection{The abstract setting}\label{sec:abstract}

We consider a PDE set in some open, connected set $\Omega \subset \mathbb{R}^d$ ($d \ge 1$). We assume that the problem can be cast in the following general abstract weak form.
Let $W$ and $V$ be two reflexive separable Banach spaces. Define a  bounded bilinear form $\aform: W \times V \mapsto \mathbb{R}$,  satisfying 
\begin{equation}\label{eq:continuity}
	\aform{(w,v)} \le M \|w\|_W \|v\|_V, \quad \forall w \in W , v \in V,
\end{equation}
and let   $\ff: V \mapsto 
\mathbb{R}$ be a bounded linear form.
We consider the abstract problem: find $u \in W$ such that  
\begin{equation}\label{eq:weakeq}
\aform(u,v) = \ff(v), \quad \forall v \in V.
\end{equation}
As in \cite{BYZZ20}, we rewrite \cref{eq:weakeq} as the following minimization problem
\begin{equation}
		\label{eq:minabstract}
		 u =
		 \underset{w \in W}{\rm \mbox{argmin}} \, 
		 \sup_{\substack{v \in V, v{ \neq 0}}}\frac{\ff(v) -\aform(w,v)}{\| v \|_V} =
		 \underset{w \in W}{\rm \mbox{argmin}} 
	\,	\|u - w\|_{op},
	\end{equation}
where we define
\begin{equation}\label{defopnorm}
\| w \|_{op} := \sup_{\substack{v \in V,v{ \neq 0}}}\frac{\aform(w,v)}{\| v \|_V}.
\end{equation}
We assume  \eqref{eq:weakeq} admits a unique solution, satisfying the following stability estimate
	\begin{equation}\label{eq:stability}
		\|u\|_W \leq C \|\ff \|_{V'}.
	\end{equation}
 This is for instance the case if the form satisfies the assumptions of the Banach-Necas-Babuska
				 theorem, or if it satisfies the more general condition of the Lions theorem, complemented by suitable compatibility conditions on $\ff$ (see \cite[Theorem 2.6 and Lemma A.40]{ern2004theory}). It is straightforward to show that, under such an assumption,  the  solution of problem \eqref{eq:weakeq} coincides with the unique minimizer of \cref{eq:minabstract}. 

\newcommand{\Wtheta}{W_{\bs\theta}}
\newcommand{\Wzeta}{V_{\bs{\eta}}}

In principle, for any function class $\Wtheta$, parametrized by a parameter set $\parameterT$,
we can approximate the solution $u$ by solving the semi-discrete problem:
\begin{equation}
	\label{eq:abstract_semidiscrete}
	\tilde u_{\bs{\theta}}^* =\underset{w_{\bs{\theta}} \in {\Wtheta}}{\rm \mbox{argmin}} \,
\|u - w_{\bs{\theta}}\|_{op}.
\end{equation}
Remark that allowing the test space $V$ to be different from the space $W$, where the solution is sought, makes the above formulation extremely flexible, allowing it to cover a wide range of situations, such as the ones where a partial differential equation, written in the form 
	\begin{equation}\label{eq:pde}
		\azang(u,w) = f(w), \qquad \forall w \in W_0 \subset W,
	\end{equation}
	($W_0$ denoting some closed subspace of $W$), is complemented by a constraint:
	\begin{equation}\label{eq:data}
		\bzang(u,\chi) = g(\chi), \qquad \forall \chi \in X,
	\end{equation}
where $X$ is a third reflexive separable Banach space.
	Such a situation falls in our abstract framework, with $V = W_0 \times X$, if we set, for $v = (w_0,\chi) \in V$,
	\[
	\aform(w,v) = a(w,w_0) + \beta b(w,\chi), \qquad \ff(v) = f(w_0) + \beta g(\chi),
	\]
	($\beta$ being a parameter weighting the constraint with respect to the equation). 		In such a case, the $\| \cdot \|_{op}$ norm satisfies
	\[
	\| w \|_{op} = \sup_{(v,\chi) \in W_0 \times X} \frac{\aform(w,(v,\chi))}{(\| v \|^2_{W_0} + \| \chi \|^2_{X})^{1/2}} \simeq  \sup_{v \in W_0} \frac{a(w,v)}{\| v \|_{W_0}} +\beta \sup_{ \chi \in  X}  \frac{ b(w,\chi)}{ \| \chi \|_{X}}.
	\]	
	Typically,  as we shall see below, \eqref{eq:data} could represent the imposition of essential boundary conditions. It could also represent some other form of constraint, such as the ones encountered in data assimilation problems subject to the heat or wave equation (see \cite{BFMO21} or \cite{BO18} where $b$ is the $L^2$-scalar product over some subset $\omega \subset \Omega$ \cite{deshmukh2016neural, moeini2012wave}). 

\subsection{The WAN method}
Let $W$ denote, throughout this section, the $H^1(\Omega)$ space with norm  $\| v \|_W$ defined as $\| v \|_W^2 = (\nabla v,\nabla v)_{\Omega} + (v,v)_\Omega$, where $(\cdot,\cdot)_\Omega$ denotes the $L^2(\Omega)$ scalar product. We consider an elliptic partial differential equation, endowed with a Dirichlet boundary condition, that we write in the form
\begin{equation}\label{zangstrong}
\Azang(u) = f, \qquad \Bzang(u) = g,
\end{equation}
where $\Azang$ is a second order partial differential operator and $\Bzang$ is the trace operator.
 Bao et al. propose in  \cite{BYZZ20}  to rewrite \cref{zangstrong} as a minimization problem in a suitable dual space. To this aim, the so-called operator norm is introduced, defined as
\begin{equation}\label{classicopnorm}
	\|\Azang(v)\|_{H^{-1}(\Omega),op}  := \sup_{\substack{\varphi \in H^1_0(\Omega)\\ \varphi \neq 0}} 
	\dfrac{\azang(v,\varphi)}{\|\varphi\|_{W}},
\end{equation} 
where $\azang: H^1(\Omega) \times H^1_0(\Omega) \mapsto \mathbb{R}$, $\azang(w,\varphi) =  (\Azang(u),\varphi)_\Omega$, is the bilinear form corresponding to the operator $\Azang$.
 We immediately see that, provided the form $\azang$ is continuous on $H^1(\Omega)\times H^1_0(\Omega)$, such norm is well defined (indeed, it coincides with the standard $H^{-1}(\Omega)$ norm).
{The idea of  \cite{BYZZ20} was then to combine the residual in such a norm with  a boundary penalization term aimed at weakly  imposing the boundary conditions (rather than enforcing them exactly),} and consider the following minimization problem: 
\begin{equation}
	\label{eq:netexa}
	u^* =\underset{w \in W}{\rm \mbox{argmin}} 
	\left(\|\Azang(u - w)\|_{H^{-1}(\Omega),op} + \beta \| g - w \|_{L^2(\partial\Omega)} \right).
\end{equation}
Setting $V = H^1_0(\Omega) \times L^2(\partial\Omega)$, and 
\[\aform(w,[v,\chi] ) = \azang(w,v) + \beta (w,\chi)_{\partial\Omega}, \qquad \ff([v,\chi]) = (f,v)_\Omega +\beta (g,\chi)_{\partial \Omega},\]
this problem can be  rewritten in the form \cref{eq:minabstract}. {At the continuous level, problem \cref{eq:netexa} is, in some sense, equivalent to  \cref{zangstrong}.  Indeed, we observe that 
the unique solution of \cref{zangstrong} annihilates both $\| A(u - w) \|_{H^{-1}(\Omega),op}$ and $\| g - w \|_{L^2(\partial\Omega)}$, implying existence. { Then, the values of the minimum is zero, and any other $H^1(\Omega)$ function minimizing the boundary penalized residual can be easily seen to be the solution to \cref{zangstrong}, thus obtaining uniqueness. }Trivially, { as it coincides with the solution of \cref{zangstrong}},  the solution of \cref{eq:netexa} satisfies $\| u^* \|_{W} \lesssim \| f \|_{H^{-1}(\Omega)} + \| g \|_{H^{1/2}(\partial\Omega)} \lesssim  \| f \|_{H^{-1}(\Omega)} + C(g) \| g \|_{L^2(\partial\Omega)}$, with $C(g) = \| g \|_{H^{1/2}(\partial\Omega)}/\| g \|_{L^2(\partial\Omega)}$, which is a stability bound of the form  \eqref{eq:stability}, though with a constant depending on $g$ {(whether such a constant is large or not depends on the frequency content of $g$: if $g$ is not oscillating, such a constant is of order one, but it can be large if $g$ presents high frequency oscillations).} However, such a formulation does not entirely fall in the abstract setting of Section \ref{sec:abstract}, since, for $V=H^1_0(\Omega) \times L^2(\partial \Omega)$, the bilinear form $\aform$ does not satisfy the boundedness assumption \cref{eq:continuity}.}
It is therefore natural to consider the following minimization problem, where the boundary penalization term is
measured in the $H^{1/2}(\partial(\Omega)) = (H^{-1/2}(\partial\Omega))'$ norm
\begin{equation}
	\label{eq:netexonehalf}
	u^* =\underset{\rm w \in W}{\rm \mbox{argmin}} 
	\left(\|\Azang(u - w)\|_{H^{-1}(\Omega),op} + \beta \| g - w \|_{H^{1/2}(\partial\Omega)} \right).
\end{equation}
It is not difficult to see that also this problem can be written in the form \cref{eq:minabstract},
this time with $V = H^1_0(\Omega) \times H^{-1/2}(\partial\Omega)$.  {Thanks to the choice of the correct norm for the boundary penalization term,} problem \cref{eq:netexonehalf} falls within our abstract framework of \cref{sec:abstract},  it is well posed, and equivalent to \cref{zangstrong}. It will serve as a starting point for the boundary condition treatment we will propose in  \cref{sec:Dirichlet}.

\begin{remark} {In the very first version of the WAN method, see \cite{ZBYZ20}, the authors actually proposed a different definition of the  operator norm, namely they  defined the dual norm involved in the minimization problem as
$
\|\Azang(v)\|_{L^{2}(\Omega),op} := \displaystyle{\sup_{\substack{\varphi \in H^1_0(\Omega), \varphi \neq 0}} 
\dfrac{\azang(v,\varphi)}{\|\varphi\|_{L^2(\Omega)}}}$. It should be noted that this norm is not generally well defined at the continuous level, and to remedy this, the different normalization in \eqref{classicopnorm} was proposed in \cite{BYZZ20}.}
 Remark that the notation used for such a norm in  \cite{ZBYZ20} was $\| \Azang(v) \|_{op}$, while we use the notation $\| \cdot \|_{op}$ with a different meaning, see \eqref{defopnorm}.
\end{remark}

{\begin{remark}
	We remark that replacing the natural norms $H^{1}(\Omega)$ and $H^{1/2}(\partial\Omega)$ in, respectively, \cref{classicopnorm} and \cref{eq:netexa} with the corresponding $L^2$-norm results in two ``variational crimes'' with fairly different features. In both cases the natural norm is replaced by a weaker norm but in the first case the replacement happens in the denominator. The resulting term $\| A(w) \|_{L^2(\Omega),op}$, $w \in W$, is not necessarily well defined (as this would  require $w \in H^2(\Omega)$). This is essentially the same residual quantity as that minimised in so-called PINN methods \cite{DP94,raissi2019physics}. In the second case, the ``variational crime'' is somewhat less severe:  all the quantities involved in the minimization problem \cref{eq:netexa} are well defined, though, as we already pointed out,  the boundedness assumption \cref{eq:continuity} does not hold. 
\end{remark}
}

In the WAN method, the discretization for 
either \cref{eq:netexa} or \cref{eq:netexonehalf} 
is performed by replacing the spaces $H^1(\Omega)$ and $H^1_0(\Omega)$ by, respectively, their discrete  counterparts ${\Wtheta} \subset H^1(\Omega)$ and $\Wzeta  \subset H^1_0(\Omega)$, where ${\Wtheta}$ and $\Wzeta $ are two fixed architecture neural network function classes, parameterized by parameter sets
$\parameterT$ and $\parameterZ$.
The discretization is carried out via a discrete operator norm, defined, for any $w \in H^1(\Omega)$, as
\begin{equation}\label{op-eta-norm}
	\|\Azang(w) \|_{H^{-1}(\Omega),op,{\bs{\eta}}} :=  \sup_{\substack{v_{{\bs{\eta}}} \in {\Wzeta}\\ \|v_{{\bs{\eta}}}\|_V \neq 0}}  
	\dfrac{\azang(w,v_{{\bs{\eta}}})}{\|v_{{\bs{\eta}}}\|_{V}} .
\end{equation}
The discrete method can then be written, for $X$ being either $L^2(\partial\Omega)$ or $H^{1/2}(\partial\Omega)$,
\begin{equation}\label{eq:netexb}
	u_{\bs{\theta}}^* =\underset{\rm w_{\bs{\theta}} \in {\Wtheta}}{\rm \mbox{argmin}} 
	\left(\|\Azang(u - w_{\bs{\theta}})\|_{H^{-1}(\Omega),op,{\bs{\eta}}} + \beta \|w_{\bs\theta} - g \|_{X} \right).
\end{equation}

Exactly 
evaluating the functional on the right-hand side is very difficult since functions in ${\Wtheta}$ and ${\Wzeta}$, may have very different geometric structures. 
In practice, the integrals are approximated using fixed sample points or a Monte Carlo integration method \cite{hong2021priori}. The optimization is then performed using a Stochastic gradient descent method, e.g. Adam, over the parameter sets $\parameterT$ and $\parameterZ$. 
{We also note that, due to the normalization in (1.13), when $w$ is close to $u$,  the maximization problem in $v_{\bs\eta}$ becomes ill-posed, resulting in increased undesirable oscillations. We will propose a possible remedy in Section \ref{sec:WAN_alternative}.

\section{Analysis of the WAN method}\label{sec:WAN_analysis}
This section will frame and analyze the WAN method in an abstract framework. We aim to provide insight into choosing
the approximation and adversarial networks to ensure the resulting method's stability and optimality. {For simplicity, we will perform the analysis based on \cref{eq:netexb} without considering the errors caused by the Monte Carlo and gradient descent methods.}

We define the WAN method in the abstract framework as follows. Letting $\Wzeta \subset V$ denote a function class
 parametrized by a parameter set $\parameterZ$, we introduce the discrete version of the $\| \cdot \|_{op}$ norm on $W$, defined as  
\begin{equation}\label{def_op_eta_abstr_norm}
	\| w \|_{op,{\bs{\eta}}} :=  \sup\limits_{\substack{v_{{\bs{\eta}}} \in {\Wzeta}\\ \|v_{{\bs{\eta}}}\|_V \neq 0}}  
\dfrac{\aform(w,v_{{\bs{\eta}}})}{\|v_{{\bs{\eta}}}\|_{V}} .
\end{equation}
We observe that, for all $w \in W$ we have that
\begin{equation}\label{eq:contnormop}
\| w \|_{op,\bs\eta} \leq \| w \|_{op} \leq M \| w \|_{W}.
\end{equation}
The fully discrete problem then reads
\begin{equation}\label{eq:neteqred}
	\wstar = \underset{w_{\bs{\theta}} \in {\Wtheta}}{\rm \mbox{argmin}} \,\,
	\|u - w_{\bs{\theta}}\|_{op,{\bs\eta} }.
\end{equation}

 In our analysis, a key role will be played by the function class of differences of elements of the approximation network $\Wtheta$:
 	\begin{equation}\label{networks-spaces}
 	S_{{\bs{\theta}}}:=\{w_{1, {\bs{\theta}}} - w_{2,{\bs{\theta}}}, \; w_{1, {\bs{\theta}}} , w_{2,{\bs{\theta}}} \in \Wtheta\}.
\end{equation} 
We will first consider the case of coercive problems and then tackle problems only known to satisfy the stability \cref{eq:stability}. 

\newcommand{\twn}{\widetilde w^n_{\bs\theta}}

\newcommand{\wn}{w^n_{\bs\theta}}

\subsection{Coercive problems} 
Let us at first consider the case $V = W$, and assume that the bilinear form $\aform$ is coercive, i.e., there exist $\alpha >0$ such that
\begin{equation}\label{eq:coercivity}
\alpha \|\phi\|_W^2 \le \aform{(\phi, \phi)}.
\end{equation} 
We make the following assumption on the networks $\Wtheta$ and $\Wzeta$:
\begin{equation}\label{assumptionZT-1}
\Wtheta \cup S_{\bs\theta} \subseteq \Wzeta.
\end{equation}
Observe that if $0 \in W_{\bs\theta}$, we have that $\Wtheta \cup S_{\bs\theta} = S_{\bs\theta}$. We start by remarking that, as the functional $w \to \| u - w \|_{op,\bs\eta}$, with $u \in W$ given, is bounded from below by $0$, we have that
\[ \minimum :=
\inf_{w_{\bs\theta} \in {\Wtheta}}	\| u - w_{\bs\theta} \|_{op,\bs\eta} \geq 0.
\]
By the definition of infimum, there exist a sequence  $\{w_{\bs\theta}^n\}$ with $w_{\bs\theta}^n \in W _{\bs\theta}$ 
such that 
\begin{equation} 
	\label{eq:minseq}
	 \lim_{n\to \infty} \| u - \wn   \|_{op,\bs\eta} =
\inf_{w_{\bs{\theta}} \in {\Wtheta}} \| u - w_{\bs\theta} \|_{op,\bs\eta}.
\end{equation}
We call a sequence satisfying \eqref{eq:minseq}  a minimizing sequence for \cref{eq:neteqred}.
We have the following lemma, where $\text{cl}^{seq}_w( {\Wtheta})\subseteq W$ denotes the weak sequential closure of ${\Wtheta}$ in $W$(see  \cite{ostrovskii2001weak}).

\begin{lem} \label{lem:existence}
	Let $\{\wn \} $ be a minimizing sequence for \cref{eq:neteqred}. Then, under assumption \eqref{assumptionZT-1}, there exists a subsequence weakly converging to an element
 $\wstar  \in \mathrm{cl}^{seq}_w( {\Wtheta})$ satisfying
		\[
	\| u - \wstar \|_{op,\bs\eta} \leq \inf_{w_{\bs\theta} \in {\Wtheta}} \| u - w_{\bs\theta} \|_{op,\bs\eta}.
	\]
\end{lem}

\begin{proof} 
Thanks to \cref{eq:coercivity,assumptionZT-1} it is not difficult to see that  the sequence $\{w_{\bs\theta}^n\}$ is bounded in $W$, and it therefore admits a weakly convergent subsequence $\{\twn \}$.
	We let $\wstar \in W$ denote the weak limit of  $\{\twn \}$. Let now $A^T: W \to W'$ be defined as
	$
	\langle A^T v, w \rangle = \aform(w,v),
	$
 with $\langle \cdot,\cdot \rangle$ denoting the duality pairing.	We have, by the definition of weak limit,
 	\begin{equation}
 		\| u - \wstar \|_{op,\bs\eta} =
 		\sup_{{v_{\bs\eta}\in {\Wzeta}}\atop{v_{\bs\eta}\neq 0}} \frac{
 			\langle A^T v_{\bs\eta} , u - \wstar \rangle
 		}{\| v_{\bs\eta} \|_{V}} 
 	= 
 		\sup_{{v_{\bs\eta}\in {\Wzeta}}\atop{v_{\bs\eta}\neq 0}} \frac{ \lim\limits_{n\to \infty}
 			\langle A^T v_{\bs\eta} , u - \twn , \rangle
 		}{\| v_{\bs\eta} \|_{V}}  .
 \end{equation}
Now, for any $v_{\bs\eta} \in {\Wzeta}$, $v_{\bs\eta} \neq 0$, we have 
\begin{equation*}
	\begin{split}
\frac{ \lim\limits_{n\to \infty}
\langle A^T v_{\bs\eta},u - \twn  \rangle }{\| v_{\bs\eta} \|_{V}}  
&\leq  \lim_{n\to \infty} \sup_{{v'_{\bs\eta}\in {\Wzeta}}\atop{v'_{\bs\eta}\neq 0}} \frac{ 
\aform (u - \twn ,  v_{\bs\eta}')  }{\| v'_{\bs\eta} \|_{V}}  =
\lim\limits_{n\to \infty} \| u - \twn  \|_{op,\bs\eta} = \minimum,
\end{split}
\end{equation*}
whence
$\| u - \wstar \|_{op,\bs\eta} \leq \sigma^*
$.
\end{proof}

We now prove Cea's lemma of best approximation for WAN on coercive problems. 

\begin{lem}\label{lem:Cea} Let assumption \eqref{assumptionZT-1} hold, and	
let $u$ be the solutions to  \cref{eq:weakeq} and  $ u_{\bs{\theta}}^* \in \mathrm{cl}^{seq}_w({\Wtheta})$ be the weak limit of a weakly convergent minimizing sequence $\{\twn \}$ for  \cref{eq:neteqred}.  
Then we have the following error bound:
\begin{equation}\label{Cea's result}
	\begin{split}
		\|u - \wstar \|_W  \le \left(1+ \frac{2 M}{\alpha}\right)
		\inf_{w_{\bs{\theta}} \in {\Wtheta}} \| u -  w_{\bs{\theta}} \|_{W}.
	\end{split}
\end{equation}
\end{lem}
\newcommand{\Riesz}{\mathfrak{R}}

\begin{proof} 
	We start by observing that \eqref{assumptionZT-1} implies that, for any two elements $w_{1,\bs\theta}$ and $w_{2,\bs\theta}$ of ${\Wtheta}$ it holds that
	\begin{equation}\label{coercivity-0}
	\alpha \| w_{1,\bs\theta} - w_{2,\bs\theta} \|_{W} \leq \sup_{{v_{\bs\eta}\in {\Wzeta}}\atop{v_{\bs\eta}\neq 0}} \frac{\aform(w_{1,\bs\theta}-w_{2,\bs\theta}, v_{\bs\eta})}{\| v_{\bs\eta} \|_W}  = \| w_{1,\bs\theta} - w_{2,\bs\theta} \|_{op,\bs\eta}.
	\end{equation}
Let  now $e^* = u - u_{\bs{\theta}}^*$, and let  $w_{\bs{\theta}}$ be an arbitrary element in $\Wtheta $. 
 Letting $(\cdot,\cdot)_W$ denote the scalar product in $W$ and $\Riesz: W \to W'$ denote the Riesz isomorphism, we have
	\begin{equation*}
		\begin{split}
	&\| \wstar - w_{\bs\theta}  \|_W =
	 \frac{(\wstar - w_{\bs\theta},
		\wstar - w_{\bs\theta}
		)_W}{\| \wstar - w_{\bs\theta} \|_W} = \frac{\langle \Riesz	(\wstar - w_{\bs\theta}),\wstar - w_{\bs\theta}
		\rangle}{\| \wstar - w_{\bs\theta} \|_W} \\
	=& \lim_{n \to \infty}  \frac{
		\langle 
		 \Riesz	(\wstar - w_{\bs\theta}),
		\twn  - w_{\bs\theta} \rangle}{
		\| \wstar - w_{\bs\theta} \|_W  
	} 
\leq \lim_{n\to\infty} \| \twn  - w_{\bs\theta} \|_{W} \leq \alpha^{-1} \lim_{n\to \infty} \| \twn  - w_{\bs\theta }\|_{op,\bs\eta}.
\end{split}
\end{equation*}
Note that we used \cref{coercivity-0} for the last bound. Adding and subtracting $u$ in the right hand side and using \cref{eq:neteqred} and \cref{eq:continuity}, we have
\begin{multline}\label{coercivity-2}	\| \wstar - w_{\bs\theta}  \|_W 	\leq \alpha^{-1} \lim_{n\to \infty} \ \| u- \twn   \|_{op,\bs\eta} + \alpha^{-1} \| u - w_{\bs\theta} \|_{op,\bs\eta}\\
	 \leq \alpha^{-1}
\inf_{w'_{\bs\theta} \in {\Wtheta}} \| u - w'_{\bs\theta} \|_{op,\bs\eta}  + \alpha^{-1} \| u - w_{\bs\theta} \|_{op,\bs\eta} \leq \frac 2 \alpha   \| u - w_{\bs\theta} \|_{op,\bs\eta}.\end{multline}
Since $w_{\bs{\theta}} \in \Wtheta $ is arbitrary, using \cref{eq:contnormop} and a triangle inequality we get  \eqref{Cea's result}.
\end{proof}

Generally, the weak solutions to  \cref{eq:netexb}, defined as the weak limits of minimizing sequences for the right hand side in ${\Wtheta}$, are not necessarily unique. 
{Moreover,  the solution of the minimization problem \eqref{eq:neteqred} itself might not lie in $W_{\bs\theta}$, but only in its weak sequential closure. In such a case, it can be proven (see \cite{petersen2021topological}) that the sequence of parameters in $\parameterT$ resulting from the minimization procedure is unbounded, which results in numerical instability.
A possible remedy (see \cite{BYZZ20}) is to restrict both maximization in $V_{\bs\eta}$ and minimization in $W_{\bs\theta}$ to subsets of $V_{\bs\eta}$ and $W_{\bs\theta}$ corresponding to parameters in $\parameterZ$ and $\parameterT$ with euclidean norm bounded by a suitable constant $B$. In such a case, 
one can apply standard calculus results to prove the existence of a minimizer $w_{\bs\theta} \in \Wtheta$. However, finding an appropriate choice of $B$ remains a challenging problem. A too-small value of $B$ will result in poor approximation regardless of the network's approximation capability, and if $B$ is very large, it ultimately serves no purpose. }
Lemma \ref{lem:Cea} does, instead, guarantee that even when multiple weak solutions exist, they all provide a quasi-best approximation of $u$ in $W$. Moreover, we can obtain a quasi-best approximation to $u$ within the approximation class $\Wtheta$ by taking entries of any minimizing sequence sufficiently close to convergence. Indeed, for any minimizing sequence $\{\twn \}$,
 given $\varepsilon > 0$ we can choose $k$ such that 
\[ \| u - \tw^k \|_{op,\bs\eta} \leq  \inf_{w_{\bs\theta} \in {\Wtheta}} \| u - w_{\bs\theta} \|_{op,\bs\eta}   + \varepsilon.\]
Then, by \cref{coercivity-2,Cea's result} we have
\[
\| u - \tw^k \|_{W} \leq \| u - \wstar \|_{W} + \| \wstar - \tw^k \|_{W}  \lesssim \inf_{w_{\bs\theta}\in {\Wtheta}} \| u - w_{\bs\theta} \|_W  + \varepsilon,
\]
meaning that any minimizing sequence does approximate the solution $u$ in the norm $\| \cdot \|_W$ within the accuracy allowed by the chosen neural network class architecture in a finite number of steps. {It is important to observe that, under proper assumptions,  the cost functional is equivalent to the $W'$ norm of the residual, thus providing a reliable a posteriori error bound. 
Moreover, by \cref{coercivity-2}, the cost functional evaluated on $w_{\bs\theta}$ provides an upper bound for the discrepancy, in $W$, between $w_{\bs\theta}$ and the weak limit $w^*_{\theta}$, and can then be leveraged to devise a stopping criterion.}

\begin{remark}\label{rem:ceadeepritz}
	Since $\Wtheta $ is a function class and not a function space, \cref{assumptionZT-1} implies that $\Wzeta $ should be a richer function class than $\Wtheta $.
	When $\aform$ is coercive and symmetric, the Deep Ritz method can 
	be interpreted as choosing, in our abstract formulation, ${\Wzeta} = u - \Wtheta$. It is not difficult to check that with such a definition of $\Wzeta$, if $\aform$ is coercive, both Lemma \ref{lem:existence} and Lemma \ref{lem:Cea} still hold.
	However, in practice, numerical evidence suggests that using a separate and more comprehensive space for ${\Wzeta}$ than $u-{\Wtheta}$ enhances both numerical efficiency (faster convergence) and accuracy. 
\end{remark}

\begin{remark} To fully exploit \cref{Cea's result}, we combine it with approximation results on neural network classes. We refer to \cite{guhring2020expressivity} for a survey of the different results available in the literature, and to the references therein.  In particular, we recall that
	when $W = H^1(\Omega)$ and  $W_{{\bs{\theta}}}$ is a function class of DNN network with ReLU activation function,
	it was shown in \cite{guhring2020error} that 
	for any function $\varphi \in H^{m}({\Omega}), m >1$ and $\Omega$ is Lipschitz,
	\begin{equation}\label{eq:approx_H2}
		\min_{\varphi_{\bs{\theta}} \in V_{{\bs{\theta}}}}\|\varphi - \varphi_{\bs{\theta}}\|_{H^{1}({\Omega})} \le C(m,d) N_{{\bs{\theta}}}^{-(m-1)/d} \|\varphi\|_{H^m(\Omega)},
	\end{equation}
	where $C(m,d)\ge 0$ is a function depends on $(m,d)$ and $ N_{{\bs{\theta}}}$ is the number of neurons in the DNN network. {Combining such a bound with the quasi-best approximation estimates allows us to deduce a priori error estimates of the WAN schemes}.
\end{remark}

{\begin{remark}
While we focused our analysis on linear problems, the WAN method can be, and is, applied also in the non linear framework. Indeed, under suitable assumptions on the operator $A$ (for instance, if $A$ is monotone and Lipschitz continuous) the existence of weakly converging minimizing sequences whose weak limit satisfies the estimate of Lemma \ref{lem:Cea} carries over to the nonlinear case. A proof in the case of monotone operators is given in the online supplementary material. Beyond that, also in cases where monotonicity does not hold, numerical results will show the effectiveness of our approach (see \cref{subsec:paraPDE} and \cref{subsec:statPDE})
	\end{remark}
}
\subsection{PDE  without coercivity}\label{sec:noncoercive}

We now drop the assumption that $V = W$,  and we assume instead that 
there exists an operator $\Riesz: V \to W'$ such that
\begin{equation}\label{inf-supWV}
\inf_{w \in W} \sup_{{v \in V}\atop{v\neq 0}} \frac{\langle \Riesz v , w \rangle}{\| w \|_{W} \| v \|_V} \geq \alpha^* > 0, \qquad \| \Riesz v \|_{W'} \leq M^* \| v \|_{V}.
\end{equation}
Remark that, as we assume that problem \eqref{eq:weakeq} is well posed, a possible choice for $\Riesz$ is $\Riesz = A^T$, but choices with better stability constants $\alpha^*$ might exist. 
Moreover assume  that ${\Wzeta} \subset V$ can be chosen 
so that we have  the discrete inf-sup condition:
\begin{equation}\label{inf-sup}
\kappa		\|w_{\bs{\theta}}\|_W \le  \sup_{{v_{\bs\eta}  \in {\Wzeta}}\atop{v_{\bs\eta}\neq 0}} \dfrac{\aform{(w_{\bs{\theta}}}, v_{\bs\eta} )}{\|v_{\bs\eta} \|_V} 
		\quad \forall w_{\bs{\theta}} \in  W_{\bs\theta} \cup S_{{\bs{\theta}}},
\end{equation} 
with $S_{\bs\theta}$ defined in \cref{networks-spaces}. It is easy to see that  \Cref{lem:existence} holds with proof unchanged also in this case, which gives us the existence of a (possibly not unique) element $\wstar \in \mathrm{cl}^{seq}_w({\Wtheta})$, weak limit of a minimizing sequence $\{\twn \}$ of elements of ${\Wtheta}$,
  satisfying
 \[
 \| u - \wstar \|_{op,\bs\eta} \leq \| u - w_{\bs\theta} \|_{op,\bs\eta} \qquad \forall w_{\bs\theta} \in {\Wtheta}.
 \]

\begin{lem}
	Let ${\Wzeta}$ be chosen in such a way that assumption \eqref{inf-sup} is satisfied
for some constant $\kappa > 0$, possibly depending on ${\Wzeta}$.
Let $u$ be the solutions to  \cref{eq:weakeq} and let $ u_{\bs{\theta}}^* \in \mathrm{cl}^{seq}_w({\Wtheta})$ be the weak limit of a weakly convergent minimizing sequence $\{\twn \}$ for  \cref{eq:neteqred}.  
Then we have the following error bound:
	\begin{equation}\label{err:noncoerc}	\| u - \wstar  \|_W 	\leq \left( 1 + 2 \frac { M^*}{\alpha_*}   \frac{M}{\kappa}  \right)  \inf_{w_{\bs\theta} \in {\Wtheta}}
		\| u - w_{\bs\theta} \|_{W}.\end{equation}
\end{lem}

\begin{proof}
Let $w_{\bs{\theta}}$ be an arbitrary element of $\Wtheta $. 
Thanks to \cref{inf-supWV,inf-sup} we can write	
	\begin{equation}\label{infsuplim}
\begin{split}	\alpha_* \| \wstar - w_{\bs\theta}  \|_W 
	&\leq \sup_{{v \in V}\atop{v\neq 0}}  \frac{\langle \Riesz	v,\wstar - w_{\bs\theta}
			\rangle}{\| v \|_V} 
		= \lim_{n \to \infty} \sup_{{v \in V}\atop{v\neq 0}}  \frac{
			\langle 
			\Riesz	v, 
			\twn  - w_{\bs\theta} \rangle}{
			\| v \|_V 
		} \\
		&\leq M^*  \lim_{n\to\infty} \| \twn  - w_{\bs\theta} \|_{W} \leq \frac{M^*} \kappa \lim_{n\to \infty} \| \twn  - w_{\bs\theta }\|_{op,\bs\eta}.
\end{split}	\end{equation}
 By the same argument used for the proof of  \Cref{lem:Cea}, we then obtain that
 \begin{equation}\label{coercivity-2b}	\| \wstar - w_{\bs\theta}  \|_W 	\leq   2 \frac{M^*}{\alpha_*}  \frac 1 \kappa   \| u - w_{\bs\theta} \|_{op,\bs\eta},\end{equation}
and, consequently, by the triangle inequality, 
	\begin{equation}	\| u - \wstar  \|_W 	\leq \left( 1 + 2 \frac { M^*}{\alpha_*}  \frac M \kappa \right)  \| u - w_{\bs\theta} \|_{W}, \end{equation}
	which, thanks to the arbitrariness of $w_{\bs{\theta}}$, gives \eqref{err:noncoerc}.
\end{proof}

Like the coercive case, we can have an almost best approximation in a finite number of steps of any weakly converging minimizing sequence $\{\twn \}$. More precisely, by the same argument as for the coercive case, for all $\varepsilon > 0$ there exists a $k$ such that
\begin{equation*}
	\| u - \tw^k \|_{W} \lesssim \inf_{w_{\bs\theta}\in {\Wtheta}} \| u - w_{\bs\theta} \|_W  + \varepsilon.
\end{equation*}

\newcommand{\normequiv}[1]{\mathcal{J}(#1)}

We conclude this section by the following observation: let $\normequiv{\cdot}$ denote any functional  on $W$ equivalent to the $\| \cdot \|_{op,\bs\eta}$ norm:
\begin{equation}\label{equivalence}
c_*  \| w \|_{op,\bs\eta} \leq \normequiv{w} \leq C^*
\| w \|_{op,\bs\eta} , \qquad \forall w \in W,
\end{equation}
and consider the problem
\begin{equation}\label{eq:netequiv}
	\wstar = \underset{w_{\bs{\theta}} \in {\Wtheta}}{\rm \mbox{argmin}} \,\,
	\normequiv{u - w_{\bs{\theta}}}.
\end{equation}
Then there exists a possibly not unique  $w^\flat_\theta \in \mathrm{cl}^{seq}_w({\Wtheta})$ such that
\[ \normequiv {u - u^\flat_{\bs\theta}} \leq \inf_{w_{\bs\theta}\in \Wtheta} \normequiv{u - w_{\bs{\theta}}}.\]
 Moreover for all $u_{\bs\theta}^\flat$ such that $u_{\bs\theta}^\flat$ is the weak limit of a minimizing sequence $\{\twn \}$ for \cref{eq:netequiv}, it holds that
\[
\| u - u_{\bs\theta}^\flat \|_{W} \leq \left(1 + 2 \frac{C^*}{c_*} \frac { M^*}{\alpha_*}   \frac{M}{\kappa} \right) \inf_{w_{\bs\theta} \in \Wtheta} \| u - w_{\bs\theta} \|_W.
\]
Indeed, by \cref{equivalence}, all minimizing sequences are bounded with respect to the $\| \cdot \|_{op,\bs\eta}$ norm, and, hence, with respect to the $\| \cdot \|_{W}$ norm. Any minimizing sequence does then weakly converge to an element $u_{\bs\theta}^\flat$. Moreover, initially proceeding as in \cref{infsuplim}, thanks to \cref{equivalence}, we have, for $w_{\bs\theta}$ arbitrary,
\begin{equation*}
	\begin{split}	{\alpha_*} \| \wstar - w_{\bs\theta}  \|_W 
 &\leq \frac{M^*} \kappa \lim_{n\to \infty} \| \twn  - w_{\bs\theta }\|_{op,\bs\eta} 
 \leq \frac {M^{*}}{ \kappa} \left(
\lim_{n\to \infty} \| \twn - u \|_{op,\bs\eta} + \| u - w_{\bs\theta }\|_{op,\bs\eta} \right)\\
&\leq \frac {M^{*}}{ \kappa c_*}  \left( \lim_{n\to \infty} \normequiv{\twn - u } + \normequiv{ u - w_{\bs\theta }} \right)  
\leq \frac {2MM^{*}C^{*}}{\kappa c_*}   
\|  u - w_{\bs\theta } \|_{W}.\end{split}
\end{equation*}

	\newcommand{\weightnorm}[1]{| \mathbf{w}(#1) |_2}

\section{Two novel stabilized loss functions}
\label{sec:WAN_alternative}

To mitigate undesirable oscillations during the optimization procedure, resulting from the inexact solution of the maximization problem involved in the definition of the operator norm, we propose two alternative definitions of the cost functional that yields the same minimum, while avoiding direct normalization.  More precisely,  we introduce two new functionals on the product space $W \times V$, such that the supremum over $v\in V$, for $w \in W$ fixed, also gives, up to a possible rescaling and translation, the operator norm of $w$, while yielding a more favorable optimization problem under discretization. 
	
\subsection{Stabilized WAN method}\label{sec:WAN_augmented} 
Define
\begin{equation}\label{eq:pert_op-a}
\begin{split}	
&\vvvert w\vvvert_{op}^2 = \sup_{v \in V} \left(\aform(w,v)-\frac{\gamma_d}{2}\|v\|_{V}^2  \right)\!, \quad 
	\vvvert w\vvvert_{op,{\bs{\eta}}}^2 = \sup_{v_{\bs{\eta}} \in V_{\bs{\eta}}} \left(\aform(w,v_{\bs{\eta}})-\frac{\gamma_d}{2}\|v_{\bs{\eta}}\|_{V}^2  \right)\!,
\end{split}
\end{equation}
where $\gamma_d >0$ is a constant, and consider the following problem:
\begin{equation}\label{eq:netabstFDreg00}
	u_{\bs{\theta}}^{\sharp} =  
	\underset{{w_{\bs{\theta}} \in W_{\bs{\theta}}}}{ \mbox{argmin}} 
	\tnorm{u - w_{\bs{\theta}}}_{op,{\bs{\eta}}}.
\end{equation}

The following lemma shows  that the norms 
defined in \eqref{eq:pert_op-a} coincide with the operator norms defined in  \cref{defopnorm} and \cref{def_op_eta_abstr_norm}, up to a constant dependent on $\gamma_{d}$.

\begin{lem}\label{norm-equivalence-a} Assume that $v_{\bs\eta} \in \Wzeta$ implies $\lambda v_{\bs\eta} \in \Wzeta$ for all $\lambda \in \mathbb{R}^+$.
Then,	for any $w \in W$, there holds
	\begin{equation} \label{tnorm-bound}
		\begin{split}
			\frac{1}{2\gamma_d} \|{w}\|_{op}^2 =
			\vvvert w\vvvert_{op}^2,
			 \qquad
			\frac{1}{2\gamma_d} \|{w}\|_{op,{\bs{\eta}}}^2 =
			\vvvert w\vvvert_{op,{\bs{\eta}}}^2.
		\end{split}
	\end{equation}
\end{lem}

\begin{remark}
	From this lemma, it seems reasonable to choose, e.g., $\gamma_d=1$. However, experiments have shown that adjusting the value of $\gamma$ can effectively control the oscillations in experiments.
\end{remark}

\begin{proof} We prove the second of the two equalities, the first can be proven by the same argument. 
	For any fixed $w \in \Wtheta$ with $w \neq 0$, and for all $\varepsilon > 0$, there exists a $\overline{\varphi}^\varepsilon_w \in V_{\bs\eta}$ (depending on $\varepsilon$) with $\| \overline{\varphi}^\varepsilon_w\|_V =1$, such that
	\[
	\aform(w,\overline{\varphi}^\varepsilon_w) \ge  (1-\varepsilon)\|w\|_{op,\bs\eta},
	\]
	which, setting $\varphi_w = \gamma_{d}^{-1} \|w\|_{op,\bs\eta}  \overline{\varphi}^\varepsilon_w$, yields
	\begin{equation}
		\begin{split}
			&\sup_{\varphi_{\bs\eta} \in \Wzeta} \left(\aform(w,\varphi_{\bs\eta})-\frac{\gamma_d}{2}\|\varphi_{\bs\eta}\|_{V}^2 \right) 
			\ge
			\aform(w,\varphi_w)-\frac{\gamma_d}{2}\|\varphi_w\|_{V}^2
			\\
			=&
			\gamma_{d}^{-1} \|w\|_{op,\bs\eta}\aform(w, \overline{\varphi}^\varepsilon_w) -\frac{1}{2\gamma_d} \|w\|_{op,\bs\eta}^2
			\ge \left(\frac 1 2-\varepsilon \right)\frac{1}{\gamma_d} \|{w}\|_{op,\bs\eta}^2.
		\end{split}
	\end{equation}
By the arbitrariness of $\varepsilon$ we then obtain that $	\vvvert w\vvvert_{op,\bs\eta}^2 \geq  \|{w}\|_{op,\bs\eta}^2/(2\gamma_d)$.
 To prove the converse inequality, using Young's equality gives
	\begin{equation*}
		\begin{split}
			\sup_{\varphi_{\bs\eta} \in \Wzeta} \left(\aform(w,\varphi_{\bs\eta})-\frac{\gamma_d}{2}\|\varphi_{\bs\eta}\|_{V}^2 \right) &\le
			\sup_{\varphi_{\bs\eta} \in \Wzeta} \left(\|w\|_{op,\bs\eta} \|\varphi_{\bs\eta}\|_V-\frac{\gamma_d}{2}\|\varphi_{\bs\eta}\|_{V}^2 \right) \le
			\dfrac{1}{ 2\gamma_d} \|w\|_{op,\bs\eta}^2.
		\end{split}
	\end{equation*}
	The first part of \eqref{tnorm-bound} can be proved in a similar way.
\end{proof}

The analysis of the minimization problem \eqref{eq:neteqred} carries then over to the minimization problem \eqref{eq:netabstFDreg00}. Then, there exists at least a minimizing sequence in $\Wtheta$ weakly converging to a limit $u_{\bs\theta}^{\sharp} \in {\rm cl}^{seq}_w(\Wtheta)$, and all weak limits of minimizing sequences satisfy either bound \eqref{Cea's result} or bound \eqref{err:noncoerc}, depending on whether $\aform$ is coercive or not.

{\begin{remark}	
When $w$ approaches the true solution, we have that 
\[
v^*(w) = \underset{v \in V}{\mbox{argmax}} \big(\aform(u - w,v)-\frac{\gamma_d}{2}\|v\|_{V}^2  \big) \to 0.\]
As a consequence, depending on how large the space $W_{\bs\theta}$ is, the problem
\[
u^*_\theta =  \underset{w_{\bs\theta} \in \Wtheta}{\mbox{argmin}} \big(\aform(u - w_{\bs\theta},v^*(w_{\bs\theta}))-\frac{\gamma_d}{2}\| v^*(w_{\bs\theta}) \|_{V}^2  \big)\]
might be close to the problem
\(u^*_\theta =   \underset{w_{\bs\theta} \in \Wtheta}{\mbox{argmin}}\, \aform(u - w_{\bs\theta},v^*(w_{\bs\theta})),
\) which, in turn, if $\| u - w_{\bs\theta}\|_W$ is small, might be ill-posed and too sensitive to the errors in evaluating $v^*(w_{\bs\theta})$. The following subsection introduces a further stabilized loss function that mitigates this issue.
\end{remark}}

\subsection{A further stabilized WAN method}\label{sec:WAN_augmented-a} 
We now define the following alternative operator norm $\tnorm{\cdot}^{+}$ as follows:
\begin{equation}\label{eq:pert_op}
	\begin{split}
		&\left(\tnorm{w}^{+}_{op}\right)^{2} := \sup_{v \in V} \left(\aform(w,v)-\frac{\gamma_d}{2}\|v\|_{V}^2 
		+ \|v\|_{V} \right),\\
		&\left(\tnorm{w}^{+}_{op,{\bs{\eta}}}\right)^{2}  := \sup_{v_{\bs{\eta}} \in V_{\bs{\eta}}} \left(\aform(w,v_{\bs{\eta}})-\frac{\gamma_d}{2}\|v_{\bs{\eta}}\|_{V}^2 
		+ \|v_{{\bs{\eta}}}\|_{V} \right),
	\end{split}
\end{equation}
where $\gamma_d >0$ is a constant, and we  define the following minimization problem:
\begin{equation}\label{eq:netabstFDreg00a}
	u_{\bs{\theta}}^{\ddag} =  
	\underset{{w_{\bs{\theta}} \in W_{\bs{\theta}}}}{ \mbox{argmin}} \,
	\tnorm{u - w_{\bs{\theta}}}_{op,{\bs{\eta}}}^{+}.
\end{equation} 

The following lemma states the relation between the norms defined in \eqref{eq:pert_op} and the operator norms.
\begin{lem}\label{norm-equivalence-b}
Assume that $v_{\bs\eta} \in \Wzeta$ implies $\lambda v_{\bs\eta} \in \Wzeta$ for all $\lambda \in \mathbb{R}^+$.	For any $w \in W$, there holds
	\begin{equation} \label{norm-bound-1}
			\frac{1}{\sqrt{2\gamma_d}} (\|{w}\|_{op}+1)=
			\tnorm{w}^{+}_{op},
			\qquad
			\frac{1}{\sqrt{2\gamma_d}}( \|{w}\|_{op,{\bs{\eta}}} + 1)
		=
			\tnorm{w}^{+}_{op,{\bs{\eta}}}.
	\end{equation}
\end{lem}

\begin{proof}  We prove the second of the two equalities, the first can be proven by the same argument.
	For any fixed $w \in W$ with $w \neq 0$, and for all $\varepsilon > 0$, 
	there exists $\overline{\varphi}^\varepsilon_w \in V_{\bs\eta}$ with $\| \overline{\varphi}^\varepsilon_w\|_V =1$ that satisfies
	\begin{equation}\label{lemma5-proof1}
		\aform(w,\overline{\varphi}^\varepsilon_w)  \ge (1-\varepsilon)\|w\|_{op,\bs\eta},
	\end{equation}
	which, setting $\varphi_w = \gamma_{d}^{-1} s\|w\|_{op,\bs\eta}  \overline{\varphi}^\varepsilon_w$, for some $s >0$ to be chosen, yields
	\begin{equation*}
		\begin{split}
			&\sup_{\varphi_{\bs\eta} \in \Wzeta} \left(\aform(w,\varphi_{\bs\eta})-\frac{\gamma_d}{2}\|\varphi_{\bs\eta}\|_{V}^2 +\|\varphi_{\bs\eta}\|_{V} \right) 
			\ge
			\aform(w,\varphi_w)-\frac{\gamma_d}{2}\|\varphi_w\|_{V}^2 +\|\varphi_{w}\|_{V}
			\\
			=&
			\gamma_{d}^{-1} s\|w\|_{op,\bs\eta}\aform(w, \overline{\varphi}^\varepsilon_w)
			-\frac{s^{2}}{2\gamma_d} \|w\|_{op,\bs\eta}^2 +\gamma_{d}^{-1} s \|w\|_{op,\bs\eta}\\
			\ge&
		(1-\varepsilon)	\frac{s}{\gamma_d} \|w\|^{2}_{op,\bs\eta}
			-\frac{s^{2}}{2\gamma_d} \|w\|_{op,\bs\eta}^2 + \dfrac{s}{\gamma_{d}} \|w\|_{op,\bs\eta}.
		\end{split}
	\end{equation*}
We can choose $s$ that maximizes the term on the right hand side. By direct computations, we have 	
\begin{equation*}
		\max_{s \in \mathbb{R}^+} \left(	(1-\varepsilon)	\frac{s}{\gamma_d} \|w\|^{2}_{op,\bs\eta}
			-\frac{s^{2}}{2\gamma_d} \|w\|_{op,\bs\eta}^2 + \dfrac{s}{\gamma_{d}} \|w\|_{op,\bs\eta}\right)
			= \frac 1 {2 \gamma_{d}} ((1-\varepsilon) \|w\|_{op,\bs\eta} + 1)^2.
	\end{equation*}
	Above we have used the fact that  
	$ \max_{s \in \mathbb{R}^+} (-as^{2} + bs) = {b^{2}}/{4a}$.
	By the arbitrariness of $\varepsilon$ we then obtain $(\tnorm{w}^{+}_{op,\bs\eta})^2 \geq (\| w \|_{op,\bs\eta} + 1)^2/(2\gamma_d)$.
The converse inequality can be proved by  once again using Young's inequality.
\end{proof}
Again, the analysis of the minimization problem \eqref{eq:neteqred}, including the existence of a weakly converging minimization sequence and a best approximation bound for all weak limits of minimizing sequences, carries over to the minimization problem \eqref{eq:netabstFDreg00}. 

\section{Imposition of Dirichlet boundary conditions}\label{sec:Dirichlet}
Herein we will adapt  to the case of WANs the technique introduced in \cite{DL20} for dealing with Dirichlet boundary conditions.
The idea is to weigh the elements of the test adversarial network by multiplying a cutoff function $\phi$ {(see \cite{sukumar2022exact} and references therein)}, so that the resulting test functions are forced to be zero on the boundary. The Dirichlet boundary condition can then be imposed on the primal network using a penalty, without violating the consistency of the equation. 
For simplicity, we assume that the problem is a symmetric second-order static elliptic PDE. We also assume that the boundary $\partial \Omega$ is smooth ($C^3$ to be precise).

We first present the ideas in the simple framework of the Deep Ritz method for the case of homogeneous boundary conditions. 
We assume the operator $\mathcal{A}$ of \eqref{eq:weakeq} to be symmetric under homogeneous Dirichlet boundary conditions. Then the continuous Ritz method may be written as
\[
u = \underset{v \in H^1_0(\Omega)}{ \mbox{argmin}}  \left(0.5 \mathcal{A}(v,v) - (f,v)_{\Omega} \right).
\]
Under the smoothness assumptions on $\partial \Omega$ we know that, provided $f$ is sufficiently smooth, $u \in H^m(\Omega)$ for some $m \ge 3$ and $\|u\|_{H^3(\Omega)} \lesssim \|f\|_{H^1(\Omega)}$.
Assuming that  $w_{\bs{\theta}} \in H^1_0(\Omega)$  for all $w_{\bs{\theta}} \in W_{\bs{\theta}}$, the Deep Ritz method takes the form
\[
u_{\bs{\theta}}^* = \underset{w_{\bs{\theta}} \in W_{\bs{\theta}}}{\mbox{argmin}}   \left(0.5 \mathcal{A}(w_{\bs{\theta}},w_{\bs{\theta}}) - f(w_{\bs{\theta}}) \right).
\]
As we already mentioned, the problem with this formulation is that it appears to be very difficult to design networks that satisfy boundary conditions by construction. Instead, typically, a penalty term of the form $\lambda \| T w_{\bs\theta} \|_{\partial\Omega}$ is added to the functional on the right hand side \cite{duan2021analysis}.
The convergence to the solution $u \in H^1_0(\Omega)$ of the continuous problem is obtained by letting $\lambda \rightarrow \infty$ and enriching the network space. In the classical numerical methods, e.g., the Finite Element Method, $\lambda$ is proportional to $h^{-s}$, with $h$ being the mesh size and $s>0$ a carefully chosen exponent. With neural network methods, it is, however, not obvious how to match the dimension of the space to the rate by which $\lambda$ grows. In general, either the accuracy or the conditioning of the nonlinear system suffers.

Our idea is to build the boundary conditions into the formulation by weighting the network functions with the level set function $\phi$, where $\phi\vert_{\partial \Omega} = 0$,  $\phi\vert_\Omega > 0$, and $\phi$ behaves as a distance function in the vicinity of $\partial \Omega$.  
 The solution we look for then takes the form $\phi w_{\bs{\theta}}$ with $w_{\bs{\theta}} \in W_{\bs{\theta}}$. The Cut Deep Ritz method reads
\[
\nu_{\bs{\theta}}^* =  \underset{w_{\bs{\theta}} \in W_{\bs{\theta}}}{\mbox{argmin}} 
\left( 0.5  \mathcal{A}(\phi w_{\bs{\theta}},\phi w_{\bs{\theta}}) - f(\phi w_{\bs{\theta}}) \right).
\]
It is straightforward to show that this is equivalent to 
\begin{equation}\label{min:45}
\begin{split}
\nu_{\bs{\theta}}^* &\!=\!  \underset{w_{\bs{\theta}} \in W_{\bs{\theta}}}{\mbox{argmin}}
\left(  \mathcal{A}(u-\phi w_{\bs{\theta}},u -\phi w_{\bs{\theta}}) \right) \!=\! 
\underset{w_{\bs{\theta}} \in W_{\bs{\theta}}}{\mbox{argmin}}\!
\left(  \!\sup_{v_{\bs\theta} \in \Wtheta}\!
\frac{\mathcal{A}(u-\phi w_{\bs{\theta}},u -\phi v_{\bs{\theta}})}{
\sqrt{\mathcal{A}(u-\phi v_{\bs{\theta}},u -\phi v_{\bs{\theta}})
}} \right).
\end{split}
\end{equation}
	
  Following \Cref{rem:ceadeepritz}, and assuming, for the sake of simplicity, that the minimization problem \eqref{min:45} has a unique solution $\nu_{\bs{\theta}}^* \in W_{\bs\theta}$ we then have that
	\[
	\| u -\phi \nu_{\bs{\theta}}^* \|_{H^1(\Omega)} \leq C \inf_{w_{\bs\theta} \in \Wtheta} \| u - \phi w_{\bs\theta} \|_{H^1(\Omega)}.
	\]

It remains to show that $\phi  w_{\bs{\theta}}$ is capable of approximating $u$ in $H^1_0(\Omega)$.
To this end let $\mathcal{O}$ be some domain such that $\Omega \subset \mathcal{O}$, where $\mathcal{O}$ is a box in $\mathbb{R}^d$ and let $\tilde u$ denote a stable extension of $u$ to $\mathcal{O}$ \cite{stein70}. We assume the following on the boundary $\partial \Omega$ and $\phi$.

\begin{assumption}\label{assumption:1}
	Let $\Omega$ be a bounded domain in $\mathbb{R}^d$.
	The boundary $\partial \Omega$ can be covered by open sets $\mathcal{O}_i, i=1, \cdots, I$, and one can introduce on every $\mathcal{O}_i$ local coordinates $\xi_1, \cdots, \xi_d$ with $\xi_d = \phi$ such that all the partial derivatives  
	$\partial \xi_i^\alpha/ \partial^\alpha x$ and  $\partial x^\alpha/ \partial^\alpha \xi$ up to order $k+1$ are bounded by some $C_0 > 0$. Moreover, $\phi$ is of class $C^{k+1}$ on $\mathcal{O}$, where $k+1 \geq 3$ is the smoothness of the domain, and
	{$|\phi| \ge M_{0}$} on 
	$\mathcal{O} \setminus \cup \mathcal{O}_i$ with some $m>0$, and in  $\cup \mathcal{O}_i$, $\phi$ is a signed distance function to $\partial \Omega$. 
\end{assumption}

We further need the following Hardy-type inequality (see \cite[Lemma 3.1]{DL20}).

\begin{lem}\label{lem:Hardy}
	We assume that the domain $\Omega$ is defined by the zero level set of the smooth function $\phi$ and that  \Cref{assumption:1} is satisfied. Then for any $v \in H^{k+1}(\mathcal{O})$ such that $v\vert_{\partial \Omega} = 0$, there holds
	\[
	\left\|v/\phi \right\|_{H^k(\Omega)} \leq C \| v\|_{H^{k+1}(\mathcal{O})}.
	\]
\end{lem}
Then, as by assumption $\|\phi\|_{W^{1,\infty}(\Omega)}<C$,  combining the quasi best approximation bound given by Lemma \ref{Cea's result} with  \cref{eq:approx_H2}
we obtain the following estimate for the Cut Deep Ritz method,
with ReLU activation function:
\[
\|u - \phi \nu_{\bs{\theta}}^*\|_{H^1(\Omega)}   \lesssim  N_{{\bs{\theta}}}^{-(m-2)/d} |u|_{H^{m}}.
\]
In particular, when $m=3$, using elliptic regularity, we can bound the $H^1(\Omega)$ norm of the error with \(
N_{{\bs{\theta}}}^{-1/d} \|f\|_{H^1(\Omega)}
\).
This shows that 
the method typically requires one order more regularity of the data than typically expected. 

\newcommand{\Xp}{{X'}}

We now introduce the Cut weak adversarial network (CutWAN) method
for problems not necessarily coercive and with non-homogeneous Dirichlet boundary conditions. We let
\begin{equation}\label{eq:pert_op_cut}
	\tnorm{w}_{op, \phi,{\bs{\eta}}} := \sup_{\varphi_{\bs{\eta}} \in V_{\bs{\eta}}} \left(\aform(w, \phi \varphi_{\bs{\eta}})-\frac{\gamma_d}{2}\|\phi \varphi_{\bs{\eta}}\|_{H^1(\Omega)}^2\right)
\end{equation}
and set
\begin{equation}\label{eq:netabstFDreg}
	u_{\bs{\theta}}^\diamond =   \underset{w_{\bs{\theta}} \in W_{\bs{\theta}}}{\mbox{argmin}}  \left( \tnorm{u - w_{{\bs{\theta}}}}_{op,\phi,{\bs{\eta}}} + \|w_{\bs{\theta}} - g\|_{H^{1/2}(\partial \Omega)}  \right).
\end{equation}
Note that the cut-off function $\phi$ only multiplies the test functions in the CutWAN network method.
The boundary condition is weakly imposed by adding a penalty term on the primal network $w_{\bs{\theta}}\vert_{\partial \Omega}$. 
It is not difficult to check that this  problem falls in the abstract formulation considered at the end of  \cref{sec:noncoercive}. Here, the role of adversarial test network is played by the product space $\phi V_{\bs\theta} \times H^{-1/2}(\partial\Omega)$, and the inf-sup condition
\eqref{inf-sup} becomes
\begin{equation}\label{inf-sup-cutWAN}
	\|w\|_{W} \lesssim \tnorm{w}_{op,\phi,{\bs{\eta}}} + \|w\|_{H^{1/2}(\partial\Omega)}, \qquad \forall w \in S_{{\bs{\theta}}}.
\end{equation}
In particular, under such an assumption, we have the following  best approximation results for the CutWAN method.
\begin{lem}\label{lem:cutwan}
	Assume that the inf-sup condition \eqref{inf-sup-cutWAN} holds. Let $u_{\bs{\theta}}^\diamond$ be the weak limit of a minimizing sequence for problem  \eqref{eq:netabstFDreg}.
Then there holds
	\begin{equation}
		\|u-u_{\bs{\theta}}^{\diamond}\|_W
		\lesssim
		\inf_{w_{\bs{\theta}} \in W_{\bs{\theta}}} \| u -  w_{\bs{\theta}} \|_{W}.
	\end{equation}
\end{lem}
Note that the CutWAN method achieves  optimal convergence rates even though the test function class is multiplied by $\phi$. So the difficulty handled by the Hardy inequality in the Cut Deep Ritz method does not appear. Indeed the difficulty of controlling the levelset weighted test function is hidden in the inf-sup assumption  \cref{inf-sup-cutWAN}. A study of this condition will be the topic of future work.

	Similarly, we can also define the following algorithm.
Define
\begin{equation}\label{eq:pert_op_cut1}
	\tnorm{w}^{+}_{op, \phi,{\bs{\eta}}} := \sup_{\varphi_{\bs{\eta}} \in V_{\bs{\eta}}} \left(\aform(w, \phi \varphi_{\bs{\eta}})-\frac{\gamma_d}{2}\|\phi \varphi_{\bs{\eta}}\|_{V}^2 
	+ \|\phi \varphi_{\bs{\eta}}\|_{V}\right),
\end{equation}
and let
\begin{equation}\label{eq:netabstFDreg1}
	u_{\bs{\theta}}^\eth =   \underset{w_{\bs{\theta}} \in W_{\bs{\theta}}}{\mbox{\rm {argmin}}}  \left( \tnorm{u - w_{{\bs{\theta}}}}^{+}_{op,\phi,{\bs{\eta}}} + \|w_{\bs{\theta}} - g\|_{H^{1/2}(\partial\Omega)} \right).
\end{equation}
We refer to the above method as the shifted CutWAN method. One can also prove the best approximation results for the shifted CutWAN method similarly as in Lemma \ref{lem:cutwan}. 

\begin{remark}\label{rem:GN}
{
For computational convenience, the $H^{1/2}(\partial\Omega)$ norm in \eqref{eq:netabstFDreg} and \eqref{eq:netabstFDreg1}  can be replaced by a suitable combination of the $L^2$ norm of the function and of its tangential derivative. Indeed, using the Gagliardo-Nirenberg inequality} we have that
		\begin{equation}\label{Gagliardo-Nirenberg}
				\|w_{\bs{\theta}} - g\|_{H^{1/2}(\partial \Omega)} 
				\leq \|w_{\bs{\theta}} - g\|^{1/2}_{L^2(\partial \Omega)}\|w_{\bs{\theta}} - g\|^{1/2}_{H^{1}(\partial \Omega)}. 
		\end{equation}
{	It is not difficult to ascertain that if we replace the $H^{1/2}$ norm in \eqref{eq:netabstFDreg} and \eqref{eq:netabstFDreg1} with the right hand side of \eqref{Gagliardo-Nirenberg}, the analysis of Section \ref{sec:WAN_analysis} holds with minor changes, resulting in an error bound of the form
	\[
	\| u - 	u_{\bs{\theta}}^\diamond \| \lesssim 	\inf_{w_{\bs{\theta}} \in W_{\bs{\theta}}} \left(\| u -  w_{\bs{\theta}} \|_{W} + \|  u -  w_{\bs{\theta}} \|^{1/2}_{L^2(\partial\Omega)}
 \|  u -  w_{\bs{\theta}} \|^{1/2}_{H^1(\partial\Omega)}\right).\]
 We observe that an analogous result would hold for the plain $L^2(\partial\Omega)$ penalization as originally proposed by \cite{ZBYZ20,BYZZ20}, if an inverse estimate were to hold for  $W_{\bs{\theta}}$, allowing to bound the $H^{1}(\partial\Omega)$ norm of the boundary residual with its $L^2(\partial\Omega)$ norm times  a constant depending on $W_{\bs{\theta}}$.   Unfortunately this is generally not true, and Lemma \ref{lem:Cea} does not hold when using such a stabilization, which is, however, computationally convenient, and which we will test extensively in the forthcoming sections. {We point out that, thanks to the combination of \eqref{Gagliardo-Nirenberg} with a Cauchy-Schwartz inequality,  simply adding $\| g - w_{\bs\theta} \|_{H^1(\partial\Omega)}$ to the $L^2(\partial\Omega)$ penalized functional yields an a posteriori error estimator.  This can be evaluated upon convergence of the optimization procedure, to check if the solution obtained with the cheaper $L^2$ penalized functional, is satisfactory. If not, it can serve, in a two stage strategy, as starting point for an additional optimization procedure relying on the more expensive functionals for which our theoretical error analysis applies.}
 }
\end{remark}

\newcommand{\bfx}{{\bs{x}}}

\section{Neural Network Structures}
\subsection{Deep Neural Network (DNN) Structure}\label{Sect:DNN}
A DNN structure is the composition of multiple linear functions and nonlinear activation functions. We will use the DNN structure for $V_{{\bs{\eta}}}$. Specifically, the first component of DNN is a linear transformation $\bs{T}^l:\mathbb{R}^{n_l}\to\mathbb{R}^{n_{l+1}}$, {$l=1, \cdots, L$}, defined as follows,
\begin{eqnarray*}
\bs{T}^l(\bfx^l) = \bs{W}^l \bfx^{l} +  \bs{b}^l,\mbox{ for }\bfx^l\in \mathbb{R}^{n_l},
\end{eqnarray*}
where $\bs{W}^l=(w_{i,j}^l)\in\mathbb{R}^{n_{l+1} \times n_l}$ and $\bs{b}^{l}\in\mathbb{R}^{n_{l+1}}$ are parameters in the DNN. The second component is an activation function $\psi:\mathbb{R}\to\mathbb{R}$ to be chosen, and typical examples of the activation functions are $\tanh$, Sigmoid, and ReLU. Application of $\psi$ to a vector $\bfx\in\mathbb{R}^n$ is defined component-wisely, i.e., $\psi(\bfx)=(\psi(x_i))$, $i=1,2,\cdots,n$. The $l$-th layer of the DNN is defined as the composition of the linear transform $\bs{T}^l$ and the nonlinear activation function $\psi$, i.e., 
\begin{equation*}
\mathcal{N}^l(\bfx^l) := \psi(\bs{T}^l(\bfx^l)),
\quad l=1, \cdots, L-1.
\end{equation*}  
For an input $\bfx\in\mathbb{R}^{n_1}$, a general $L$-layer DNN is defined  as follows,
\begin{equation}\label{eqn:DNN}
\mathcal{NN}(\bfx; {\bs{\theta}}) := \bs{T}^L \circ \mathcal{N}^{L-1} \circ \cdots \circ \mathcal{N}^{2} \circ \mathcal{N}^1(\bfx),
\end{equation}
where ${\bs{\theta}} \in \mathbb{R}^N$ stands for all the parameters in the DNN, i.e.,
$
{\bs{\theta}} = \{\bs{W}^l, \bs{b}^l\}_{l=1}^{L}.
$
For a fully connected DNN, the number of parameters corresponding to ${\bs{\theta}}$ is  $N_{{\bs{\theta}}} :=  \sum_{l=1}^L n_{l+1}(n_{l} + 1)$.
We will refer $\mathcal{N}^{1}$ as the input layer, $\mathcal{N}^{i}, 1<i < L$ as the hidden layers, and $T^{L}$ as the output layer. We assume that every DNN neural network has an input, an output, and at least one hidden layer. Note that for the outer layer, there is no followed activation function.
\Cref{fig:u-DNN} shows an example of a DNN model with $5$ hidden layers with $[n_{1}, n_{2}, \cdots, n_{7}] = [6, 20, 10, 10, 10, 20,1]$.

\begin{figure}[!ht]
\centering
\includegraphics[width=0.8\textwidth]{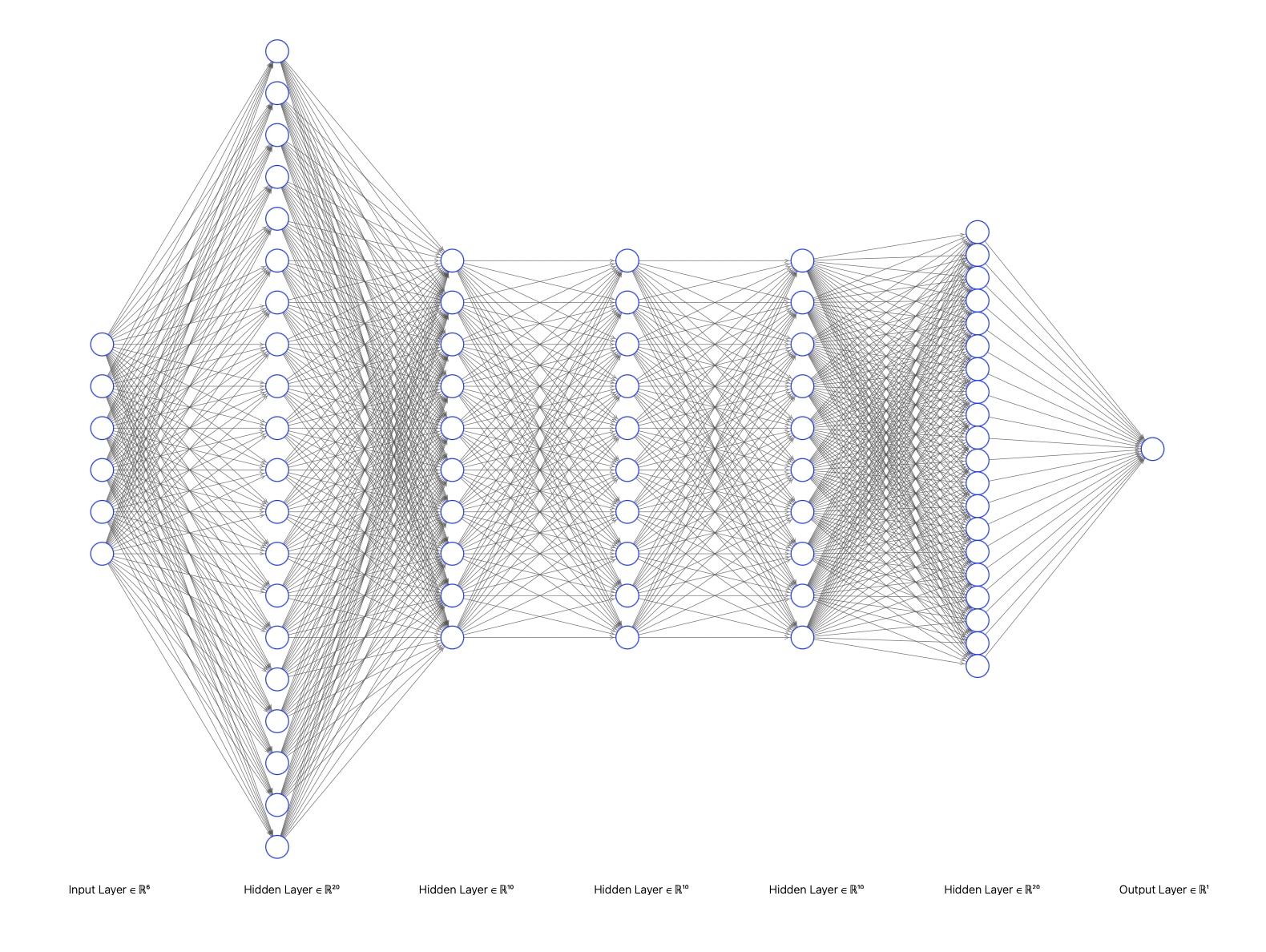} 
\caption{ A DNN network structure with $5$ hidden layers} \label{fig:u-DNN}
\end{figure} 

\subsubsection{The recursive DNN model}
In the case of a DNN model with consecutive hidden layers having an equal number of neurons, the weights and biases for those hidden layers can be easily shared due to the same data structure. We define the recursive DNN model as  DNN models that share the parameters for all consecutive hidden layers with the same number of neurons. Therefore, A recursive DNN model could have significantly fewer total parameters than the corresponding non-recursive DNN model. For instance, the non-recursive DNN model described in  \Cref{fig:u-DNN} has $931$ parameters, while its corresponding recursive model has only $711$ parameters. The contrast will become more pronounced when the number of hidden layers and hidden neurons increases. 
Our numerical results show that a recursive DNN model can benefit PDE solving. 

\subsubsection{Comments about DNN}
Although DNNs have been widely used as the primary neural network for solving PDE problems, their performance often falls short of expectations. When using DNNs within the Physically Informed Neural Networks  (PINN) and Deep Ritz methods, achieving the desired accuracy typically requires thousands of iterations due to oscillations and stagnation. The method of WAN helps the algorithm escape local minima. However, despite this improvement, the number of iterations remains in the range of several thousand, as reported in \cite{ZBYZ20} and demonstrated in our numerical results in  \Cref{sec:num}. To enhance convergence, we explore different neural structures that approximate the trial functions with more efficacy.

\subsection{XNODE model for parabolic PDE.}
It has been demonstrated in \cite{oliva2021fast} that for time-dependent parabolic problems, the XNODE model achieves much faster convergence than traditional deep neural networks. We believe this rapid convergence is attributed to the structure of the XNODE model, which emulates the residual network, and the direct embedding of the initial condition in the model. 

Consider the following parabolic PDE defined on an arbitrary bounded domain $\mathcal{D} \subset [0, T] \times \mathbb{R}^{d}$, possibly representing a time dependent spatial domain,
\begin{align}\label{pde}
\begin{cases}
\partial_t u
- \nabla \cdot A(t,\bs{x}) \nabla u + \bs{b}(t,\bs{x})\nabla u 
+c(u, \bs{x})u-f(\bs{x})=0& \text{for } (t, \bs{x}) \in \mathcal{D},\\
u(t, \bs{x})= g(t,\bs{x}) & \text{on }\partial \mathcal{D},\\
u(0,\bs{x})-h(\bs{x})=0 & \text{on }\Omega(0).
\end{cases}
\end{align}
where  $A = \{a_{ij} \}$, $\bs{b} = \{b_1, b_{2}, \cdots, b_{n} \}$, $f:\mathcal{D}\to\mathbb{R}$, $c: \mathbb{R} \times \mathcal{D} \to \mathbb{R}$ and $h: \Omega(0)  \to \mathbb{R}$ are given, with $\Omega(t):= \{\bs{x} | (t, \bs{x}) \in \mathcal{D}\}$ denote the spatial domain of $\mathcal{D}$ when restricting time to be $t$. Note that $c$ can be a non-linear function with respect to the first argument.

We now briefly introduce the XNODE model in \cite{oliva2021fast}. For simplicity, we consider a time-independent domain in this paper, i.e., 
$\mathcal{D} = [0,T] \times \Omega$, where $\Omega \subset \mathbb{R}^d$ is bounded. 

The XNODE model maps an arbitrary input $\bs{x} \in \mathbb{R}^d$ to the output $o_{\bs{x}}(t)_{t \in [0, T]} \in \mathcal{C}([0, T]$ by solving the following ODE problem: 
\begin{equation}\label{XNODE}
\begin{cases}
\frac{\mathrm{d}\bs{h}(t)}{\mathrm{d}t} = \mathcal{N}^{\text{vec}}_{{\bs{\theta}}_2}(\bs{h}(t), t, \bs{x}), \quad \bs{h}(0) = \mathcal{N}^{\text{init}}_{{\bs{\theta}}_{1}}(h(\bs{x})) \in \mathbb{R}^h.& \\
o_{\bs{x}}(t)= \mathcal{L}_{{{\bs{\theta}}}_{3}}(\bs{h}(t)). &
\end{cases}
\end{equation}
where $\mathcal{N}^{\text{vec}}_{{\bs{\theta}}_2}$ and $\mathcal{N}^{\text{vec}}_{{\bs{\theta}}_1}$  are DNN neural networks fully parameterized by $\mathcal{P}_{\bs{\theta}_{2}}$ and $\mathcal{P}_{\bs{\theta}_{1}}$  for the vector fields and the initial condition $\bs{h}(0)$ respectively. $\mathcal{L}_{{{\bs{\theta}}}_{3}}$ is a single linear layer parameterized by $\mathcal{P}_{\bs{\theta}_{3}}$. By $\Theta = ({\bs{\theta}}_1, {\bs{\theta}}_2, {{\bs{\theta}}}_{3})$ we denote the set of all trainable model parameters of the proposed XNODE model. 
Finally define 
\begin{equation}\label{xnode-output}
u_{{\bs{\Theta}}}(t, \bs{x}) := o_{\bs{x}}(t) \approx u(t, \bs{x})\quad \forall \bs{x} \in \Omega.
\end{equation}

\subsection{Pseudo-time XNODE model for static PDEs}
In this subsection, we expand the XNODE model  to handle stationary PDE problems. To simplify matters, we will focus on the following form of stationary PDE problem.
\begin{align}\label{pde-stat}
\begin{cases}
- \nabla \cdot A(\bs{x}) \nabla u({x}) + \bs{b}(\bs{x})\cdot\nabla u(\bs{x})
+c(u, \bs{x}) u-f(\bs{x})=0& \bs{x} \in \Omega = [0,1]^{d},\\
u(\bs{x})= g(\bs{x}) & \text{on }\partial  \Omega.
\end{cases}
\end{align}
The idea is to introduce a pseudo-time variable, which we choose from one of the spatial variables, $x_{i}$, to compensate for the absence of $t$, i.e., we let $t = x_{i}$ for some prefixed $i$. For simplicity, we choose $i=1$ without loss of generality. The remaining variables ${x_{i}, i=2, \cdots, d}$ will form the spatial variables in the XNODE model.
More precisely, the spatial input point for the pseudo-time XNODE model should be modified as  $\tilde{\bs{x}}= \{x_{2}, \cdots, x_{d} \}$.
Similar to \cref{xnode-output}, we now define
\begin{equation}\label{xnode-output1}
u_{{\bs{\Theta}}}(\bs{x}) = u_{{\bs{\theta}}}(x_{1}, \tilde{\bs{x}}) :=o_{\tilde{\bs{x}}}(x_{1})   \approx u(x_{1}, \tilde{\bs{x}}),
\end{equation}
where $o_{\tilde{\bs{x}}}(x_{1})$ is the numerical solution of \eqref{XNODE}.

\subsection{Loss functions}
We first recall the classical WAN loss function used in \cite{oliva2021fast}:
\begin{equation}\label{WAN-loss}
L_{\text{wan}}({\bs{\theta}}, {\bs{\eta}}) := \log\left(\frac{|(\aform(u_{\bs{\theta}})-f, \phi v_{\bs{\eta}})|^{2}}
{\| \phi v_{\bs{\eta}}\|^{2}_{L^2(\mathcal{D})}}\right)+ \alpha L^{2}_\text{init}({\bs{\theta}}) +  \beta L^{2}_\text{bdry}({\bs{\theta}}),
\end{equation}
where $\alpha$, $\gamma$ are hyperparameters as penalty terms and
\begin{equation*}
L_\text{init}({\bs{\theta}})= \| u_{\bs{\theta}}(0,\bs{x})-h(\bs{x})\|_{L^2(\Omega)}, \quad
L_\text{bdry}({\bs{\theta}})= {\| u_{\bs{\theta}}(t,\bs{x})-g(t,\bs{x})\|_{L^2([0,T] \times \partial \Omega)},}
\end{equation*}
and $\phi(\bs{x})|_{\partial \Omega} =0$. Here $u_{\bs{\theta}} \in W_{{\bs{\theta}}}$ and $v_{{\bs{\eta}}} \in V_{{\bs{\eta}}}$ where $W_{{\bs{\theta}}}$ and  $V_{{\bs{\eta}}}$ are neural network function classes parameterized by ${\bs{\theta}}$ and ${\bs{\eta}}$, respectively. In this paper, we use the classical DNN function class for 
$V_{{\bs{\eta}}}$. For $W_{{\bs{\theta}}}$, we will utilize and compare different neural network structures, which will be specified in each experiment.
When the PDE problem is static, $\alpha$ is set to $0$.

We also define the loss functions for the respective Cut-WAN and shifted Cut-WAN methods,
\begin{align}\label{AL-MAL}
\begin{split}
L_{\text{cwan}}({\bs{\theta}}, {\bs{\eta}})& = |(\aform(u_{\bs{\theta}})-f, \phi v_{\bs{\eta}})| - \gamma_{d} \|\phi v_{{\bs{\eta}}}\|_{V}^{2} +
\alpha L^{2}_\text{init}({\bs{\theta}})  +  \beta L^{2}_\text{bdry}({\bs{\theta}}),\\
L_{\text{scwan}}({\bs{\theta}}, {\bs{\eta}})& = 
L_{\text{cwan}}({\bs{\theta}}, {\bs{\eta}}) + \|\phi v_{{\bs{\eta}}}\|_{H_{0}^{1}(\Omega)}.
\end{split}
\end{align}
During computation, the integrals are estimated using the Monte-Carlo sampling.

\section{Numerical Results}\label{sec:num}
The {authors} carried out the numerical results on a personal CPU device (Apple M1 Max chip with 32 GB memory and 10 total cores). {The Adam optimization method is used for all presented numerical experiments.}
\subsection{Parabolic equations}\label{subsec:paraPDE}

\begin{ex}\label{ex1}
Following the numerical example in \cite{ZBYZ20, oliva2021fast}, we consider the following non-linear PDE problem in the form of a $d$-dimensional nonlinear diffusion-reaction equation (Eq. \eqref{eqn:problem}) defined on a bounded domain $\mathcal{D} \subset  [0,1]\times [-1, 1]^{d}$:
\begin{equation}\label{eqn:problem}
\begin{cases}
\partial_t u-\bigtriangleup u - u^2 -f=0 & \text{for }(t,\bs{x})\in \mathcal{D}\\
u-g=0 & \text{on }\partial\mathcal{D}\\
u(0,\bs{x})-h(\bs{x})=0 & \text{on }\Omega(0),
\end{cases}
\end{equation}
where the  exact solution is given by 
\begin{equation}\label{eqn:eg1_sol}
u(t,\bs{x}) = 2\sin\left(\frac{\pi}{2}x_1\right) \cos\left(\frac{\pi}{2}x_{2}\right)e^{-t}.
\end{equation}
\end{ex}


 
 The hyperparameters for the XNODE model for $u$ and DNN model for $v$ used in these experiments are listed in Table \ref{tab:1}, and their meanings are explained in \cref{appen:A}. The same hyperparameters were maintained across all experiments in Example \ref{ex1} for the XNODE model. The recursive (nonrecursive) XNODE model $u_{\bs{\theta}}$ has $1501\,(2161)$ trainable parameters, while the recursive (non-recursive) model of $V_{\bs{\eta}}$ has $5902\, (23351)$ trainable parameters.
\begin{table}[!ht]
\centering
{\small
\begin{tabular}{|l|l|l|l|l|l|l|l|}
\hline
 $d$ &$N_r$  & $N_b$ &$n_{T}$  & $K_{u}$  & $K_{\phi}$  & $\alpha$  & $\beta$  \\ \hline
 $5$  &  $4000$ & $4000$    & $20$          &  $2$            &  $1$                & $10^{7}$ & $10^{5}$   \\ \hline
 $\epsilon$ & $l_{\theta}$  & $l_{\eta}$  & $u_{\text{layers}}$  & $u_{\text{hid-dim}1}$  &$u_{\text{hid-dim}2}$  	 
  &$v_{\text{layers}}$  
  & $v_{\text{hid-dim}}$  \\ \hline
 $10^{-2}$&  $.015$& $.04$  & $8$  & $20$  & $10$  & $9$  &  $50$ \\ \hline
\end{tabular}
\caption{Hyper parameter setting for Example \ref{ex1}}
\label{tab:1}
}
\end{table}

From \Cref{tab:1}, large penalty constants for $\alpha$ and $\beta$ are utilized. We hypothesize that larger penalty constants can help strongly enforce initial and boundary conditions, which is beneficial in PDE solving using neural networks.

When utilizing the XNODE model to compute $u_{\bs{\theta}}$, the training process ceases either when the relative training error drops below $1\%$ or after a maximum of $300$ iterations. Conversely, if the DNN model is used to compute $u_{\bs{\theta}}$, the maximum number of iterations is set to $3000$.

For a comparison, we first train the models using the PINN type loss function defined as follows: 
\begin{equation}\label{PINN}
L_{\text{pinn}}(\theta) = \|\partial_t u_{\bs{\theta}}-\bigtriangleup u_{\bs{\theta}} - u_{\bs{\theta}}^2 -f\|_{L^{2}(\Omega)}+ \alpha L_\text{init}(\bs{\theta})  +  \beta L_\text{bdry}(\bs{\theta}).
\end{equation}
The $L_{\text{pinn}}$ type loss function was initially introduced in the physics-informed neural network by Raissi et al. (2019) \cite{raissi2019physics}. We have conducted experiments on $L_{\text{pinn}}$ with the random initialization, and the results are displayed in \Cref{Fig:0}. We utilized the XNODE and DNN models for both the recursive and non-recursive versions.
We note that the PINN loss function requires computing higher-order derivatives, which poses potential challenges. Firstly, the loss function becomes invalid when there is no strong solution, and secondly, computing these derivatives increases the computational time.
\begin{figure}[!ht]
    \centering
        \subfloat[XNODE(R) +$L_{\text{pinn}}$]{\includegraphics[width=0.48\textwidth]
    {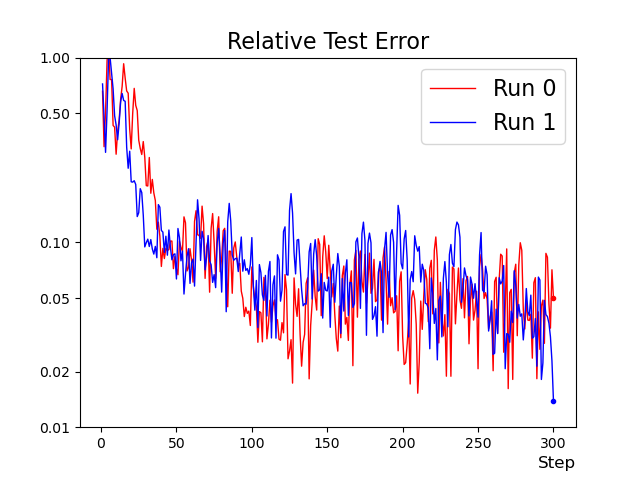} \label{fig0:ex1-PINN-XNODER}}
        \subfloat[XNODE +$L_{\text{pinn}}$]{\includegraphics[width=0.48\textwidth]{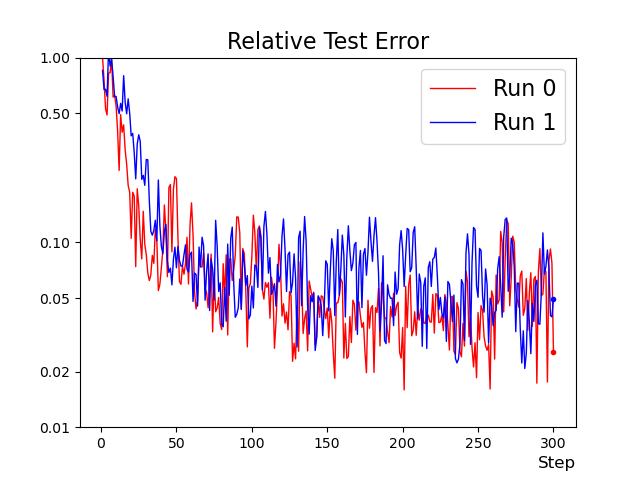}
         \label{fig0:ex1-PINN-XNODE}}\\
              \subfloat[DNN(R)+$L_{\text{pinn}}$]{\includegraphics[width=0.48\textwidth]{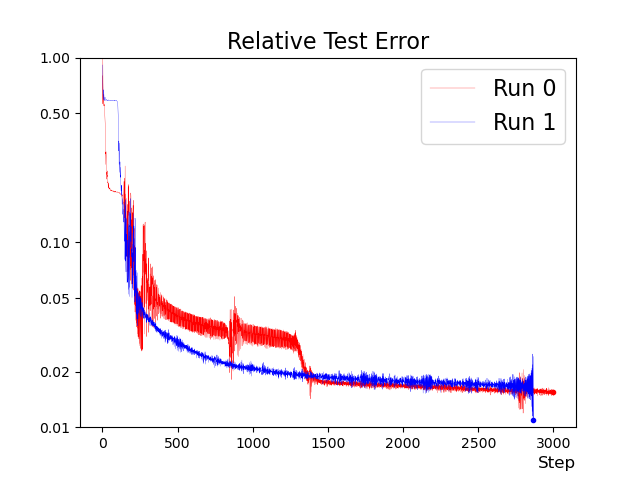}
         \label{fig0:ex1-DNN-PINNR}}
        \subfloat[DNN+$L_{\text{pinn}}$]{\includegraphics[width=0.48\textwidth]
        {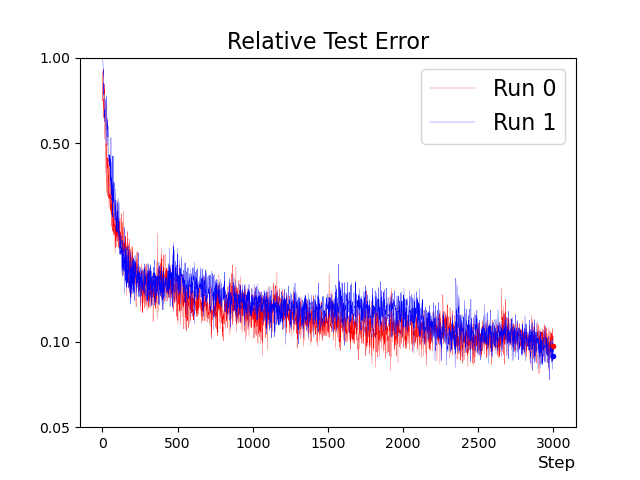}
         \label{fig0:ex1-DNN-PINN}}
\caption{ Example \ref{ex1}: Relative  $L^{2}$ Error versus Step for models using $L_{\text{pinn}}$} \label{Fig:0}
\end{figure} 
In each step, the relative error in \Cref{Fig:0} and subsequent figures is calculated using a randomly chosen test set, denoted as $X_{test}$, that is the same size as the training data sets. More precisely,
the relative error is computed as
\[
	\sum_{\bs{x_{i}} \in X_{test}} \dfrac{\sum_{i} (u(\bs{x_{i}}) - u_{\theta}(\bs{x_{i}}))^{2} }{\sum_{i} u(\bs{x_{i}})^{2}}
\]
 It's worth noting that the test set is separate from the training set but has the same size.

In \Cref{Fig:0}, 
the training time for the DNN and XNODE models is about $2$ and $8$ seconds per step. The ``(R)'' after the model denotes the recursive model.
In terms of $L_{\text{pinn}}$, it is evident that the DNN model exhibits slower convergence than the XNODE model, whether in recursive or non-recursive scenarios. When we compare figures in the right column from the left column, it is apparent that the recursive DNN model produces comparable results.

%

When using the XNODE model, from Figure \ref{fig0:ex1-PINN-XNODER} and Figure \ref{fig0:ex1-PINN-XNODE}, in both cases, relative errors approached the $7\%$ threshold within the first $50$ iterations. However, the errors then oscillate with large amplitude, requiring various steps to achieve the next level of accuracy.

In \cref{Fig:00}, we then train the models using $L_{\text{wan}}$ using the same models as in \Cref{Fig:0}.
\begin{figure}[!ht]
    \centering
        \subfloat[XNODE(R) +$L_{\text{wan}}$]{\includegraphics[width=0.48\textwidth]{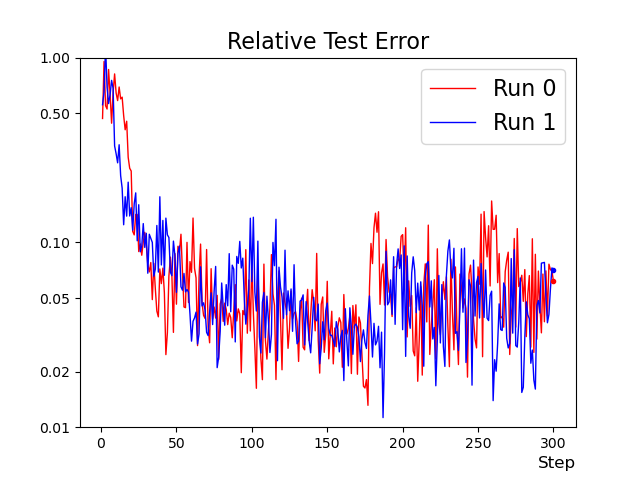}
         \label{fig00:ex1-XNODER-WAN}}
          \subfloat[XNODE+$L_{\text{wan}}$]{\includegraphics[width=0.48\textwidth]{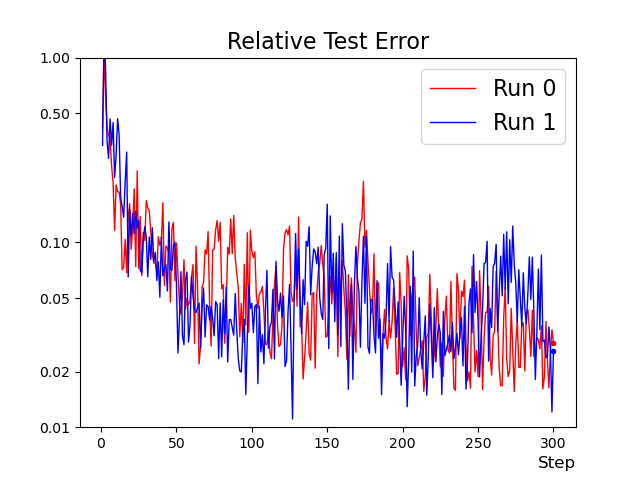}
         \label{fig00:ex1-XNODE-WAN}}\\
        \subfloat[DNN(R) +$L_{\text{wan}}$]{\includegraphics[width=0.48\textwidth]
        {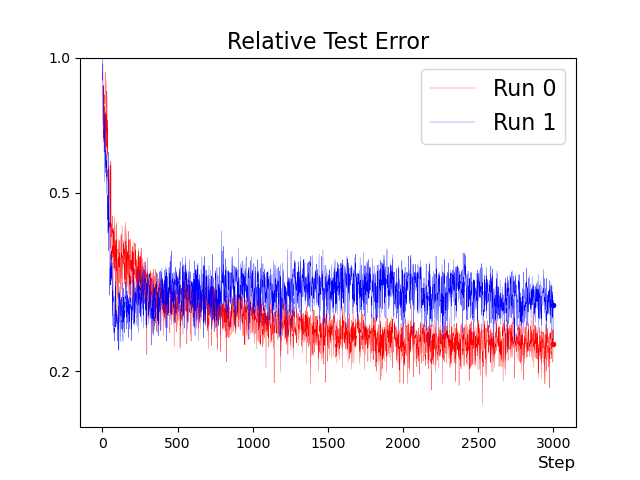}
         \label{fig00:ex1-DNNR-WANR}}
          \subfloat[DNN +$L_{\text{wan}}$]{\includegraphics[width=0.48\textwidth]
        {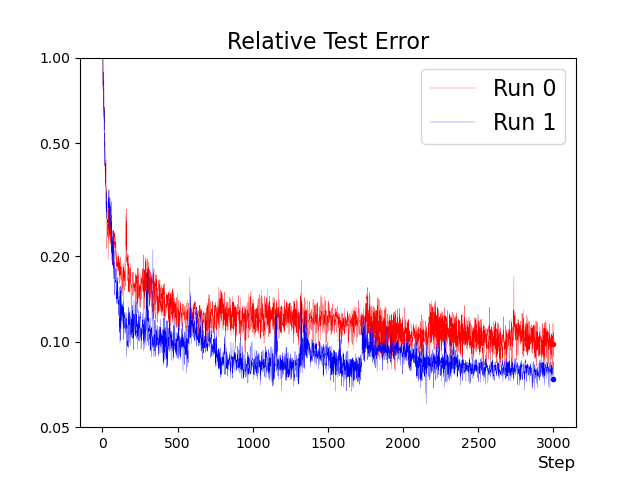}
         \label{fig00:ex1-DNN-WAN}}
\caption{ Example \ref{ex1}: Relative $L^{2}$ Error versus Step for models using $L_{\text{wan}}$ }
 \label{Fig:00}
\end{figure} 
The training time for the DNN and XNODE models is about $2$ and $6$ seconds per step.
After analyzing both \Cref{fig00:ex1-DNNR-WANR} and \Cref{fig00:ex1-DNN-WAN}, it is apparent that the utilization of 
$(L_{\text{wan}}+$ DNN) produces less desirable results compared to ($L_{\text{pinn}}+$ DNN) based on \Cref{Fig:0}. However, the combination of ($L_{\text{wan}} +$ XNODE) produces comparable results with 
($L_{\text{pinn}} +$ XNODE). This indicates that XNODE is less sensitive to the chosen objective function.

Based on the observations from  \Cref{Fig:0} and \Cref{Fig:00}, we can deduce that the utilization of the XNODE network for $u_{\bs{\theta}}$ outperforms the DNN network in both the $L_{\text{wan}}$ and $L_{\text{pinn}}$ scenarios. However, it is noteworthy that when employing the XNODE network, the loss function during the training exhibits significant oscillation after reaching a certain level of accuracy for both $L_{\text{wan}}$ and $L_{\text{pinn}}$. Additionally, in every scenario presented, the recursive model delivers comparable 
outcomes to its non-recursive counterpart.

We now evaluate the XNODE network using the loss functions $L_{\text{cwan}}$ and $L_{\text{scwan}}$ defined in \cref{AL-MAL}, with the same models as shown in \Cref{Fig:0}.
In each sub-figure in \Cref{Fig:1}, we present the results of five {out of six} consecutive and randomly initialized experiments under the specific setting to show generality. Each training step takes about $6$ seconds.

Overall, after comparing \Cref{Fig:0} and \Cref{Fig:00} with \Cref{Fig:1}, we have noticed that $L_{\text{cwan}}$ and $L_{\text{scwan}}$ show uniformly faster and {numerically more stable, i.e., less oscillations,} convergence than $L_{\text{wan}}$. Moreover, we observe {consistent/robust performance regardless of random initialization.} 
In almost all experiments, the training relative error reaches the $1\%$ relative error all within $200$ steps. 

When we compare the data in the right column to that of the left column, we notice that the recursive model performs just as well, if not better. Specifically, in the experiments depicted in \Cref{fig:ex1-AL-gamma-0.001}, the stopping criteria were met at an average of $144$ steps, with individual results of $137, 85, 177, 132$, and $190$ for experiments 0 to 4, respectively. On the other hand, the non-recursive counterpart met the stopping criteria at an average of $164$ steps, with individual results of $168, 125, 210, 153$, and $162$ for experiments 0 to 4, respectively, as shown in \Cref{fig:ex1-NRAL-gamma-0.001}.
We also note that for the $L_{\text{scwan}}$, the {comparison} results for $\gamma_{d}=0.5$ and $\gamma_{d}=0.001$ are similar in this example. {However, with $\gamma_{d}=0.001$, we notice slightly more oscillaltions than the case of  $\gamma_{d}=0.5$.}


In summary, the utilization of the XNODE network and the cutWAN and shifted CutWAN loss functions, i.e., $L_{\text{scwan}}$ and $L_{\text{cscwan}}$in \cref{AL-MAL},  has demonstrated a highly competitive model for solving high-dimensional parabolic PDE problem.
In particular, solving the $5$ dimensional non-linear parabolic problem in \cref{eqn:problem} takes only about $15$ minutes for the training to reach the $1\%$ relative error on a personal computer. Furthermore, the recursive model necessitates fewer parameters in contrast to non-recursive models. In comparison to classical numerical techniques such as the finite element method, which grows exponentially in the number of unknowns as the dimension expands, the potential benefit of our approach becomes more prominent as the disk space on a personal computer can rapidly become restricted with the classical approach.

\begin{figure}[!ht]
    \centering
        \subfloat[XNODE(R) + $L_{\text{cwan}}$ ($\gamma_{d} = 10^{-3}$)]{\includegraphics[width=0.48\textwidth]
        {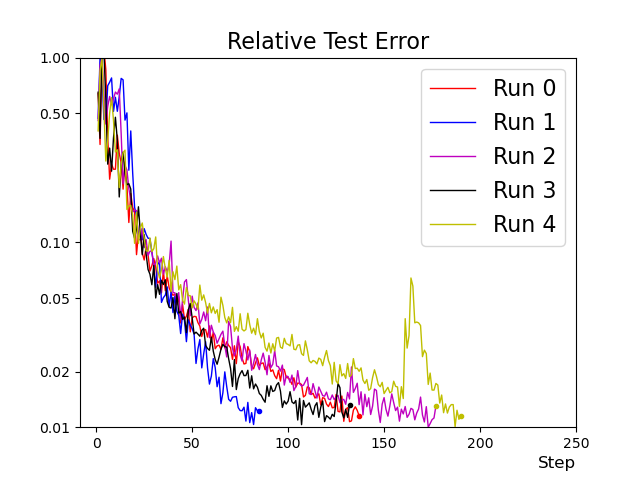}
         \label{fig:ex1-AL-gamma-0.001}}
       \subfloat[XNODE + $L_{\text{cwan}}$ ($\gamma_{d} = 10^{-3}$)]{\includegraphics[width=0.48\textwidth]
       {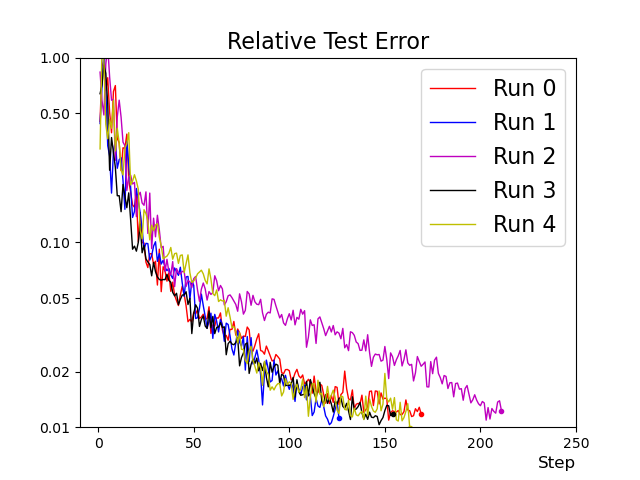}
         \label{fig:ex1-NRAL-gamma-0.001}}\\   
              \subfloat[XNODE(R) +$L_{\text{scwan}} (\gamma_{d}=0.5)$]{\includegraphics[width=0.48\textwidth]
              {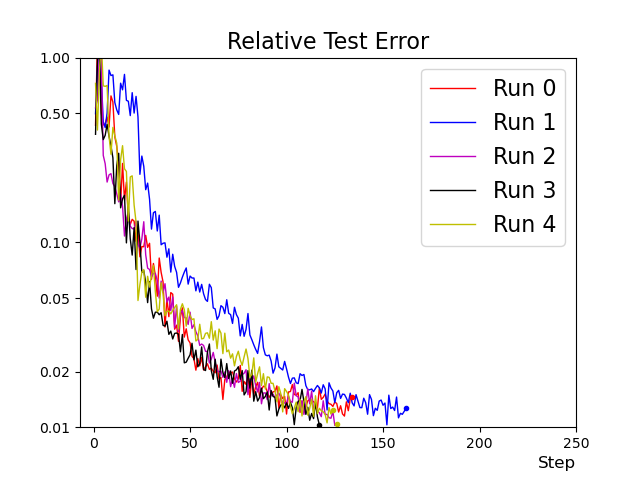}
         \label{fig:ex1-RMAL-0.5}}
                   \subfloat[XNODE +$L_{\text{scwan}} (\gamma_{d}=0.5)$]{\includegraphics[width=0.48\textwidth]
                   {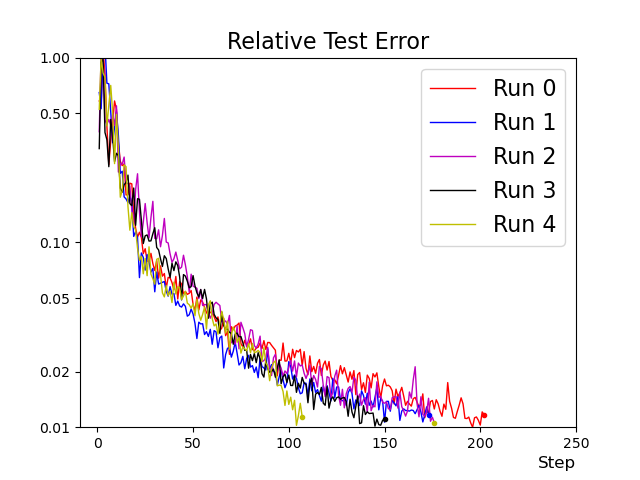}
         \label{fig:ex1-MAL-0.5}}\\
          \subfloat[XNODE(R) +$L_{\text{scwan}} (\gamma_{d}=10^{-3})$]{\includegraphics[width=0.48\textwidth]
          {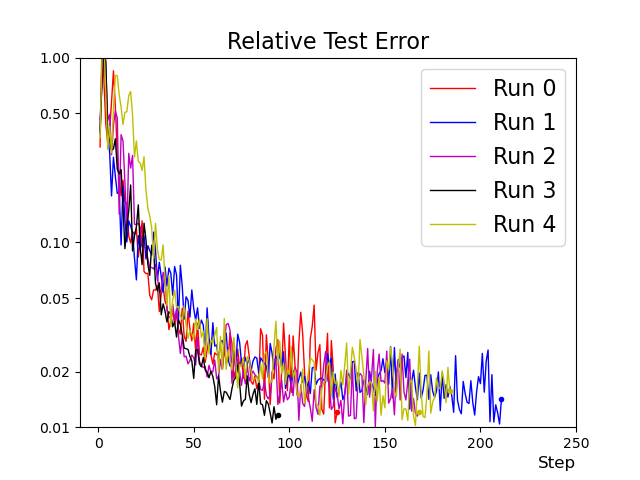}
          
         \label{fig:ex1-RMAL}}
                   \subfloat[XNODE +$L_{\text{scwan}} (\gamma_{d}=10^{-3})$]{\includegraphics[width=0.48\textwidth]
                   {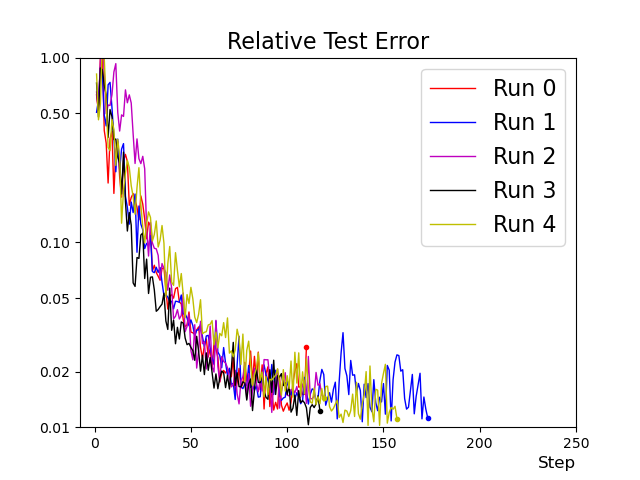}
         \label{fig:ex1-MAL}}\\
        
\caption{Example \ref{ex1}: Relative $L^{2}$ Error for XNODE models on $L_{\text{cwan}}$ and $L_{\text{scwan}}$} \label{Fig:1}
\end{figure}


{
We now consider the effect using the $H^{1/2}$ norm on the boundary based on \cref{Gagliardo-Nirenberg}.
Define
\[
\tilde L_\text{bdry}({\bs{\theta}})= \|u\|^{1/2}_{L^{2}(0,T,\partial \Omega)} 
 \|\nabla_{\Gamma} u\|^{1/2}_{L^{2}(0,T,\partial \Omega)}
\]
where $\nabla_{\Gamma} u = \nabla_{\bs{x}} u - (\nabla_{\bs{x}} u \cdot \bs{n}) \bs{n}$ is the tangential gradient of $u$ and $\bs{n}$ is the unit outer normal of $\Omega$. We test the results replacing $L_{bdry}$ in \cref{WAN-loss} by 
$\tilde L_\text{bdry}$ using the loss function $L_{\text{scwan}}$ with $\gamma_{d} = 1/2$ (see \cref{Fig:4_revision}). 
The results in \cref{fig:ex4a-RMAL-0.5} (average iteration number $=151$) and  \cref{fig:exea-MAL-0.5} (average iteration number $=173$)  are comparable to  \cref{fig:ex1-RMAL-0.5}  (average iteration number $=132$)   and \cref{fig:ex1-MAL-0.5}  (average iteration number $=161$) . 
However, using $\tilde L_\text{bdry}$ resulted in an additional duration of approximately 1 second per iteration. 
For simplicity, we will use  $L_\text{bdry}$ for future experiments. 
It is worth noting that one can use  $L_\text{bdry}$ for the former iterations and switch to the more accurate 
$\tilde L_\text{bdry}$ for better accuracy and time efficiency. This can be necessary when $g$ is of high frequency.
}

\begin{figure}[!ht]
    \centering
              \subfloat[XNODE(R) +$L_{\text{scwan}} (\gamma_{d}=0.5)$]{\includegraphics[width=0.48\textwidth]
              {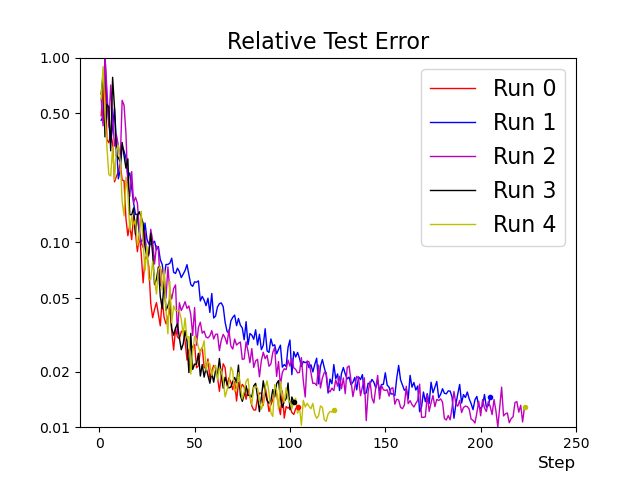}
         \label{fig:ex4a-RMAL-0.5}}
                   \subfloat[XNODE +$L_{\text{scwan}} (\gamma_{d}=0.5)$]{\includegraphics[width=0.48\textwidth]
                   {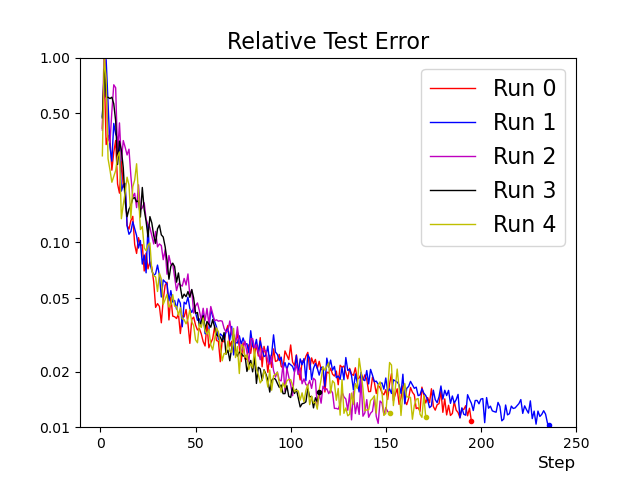}
         \label{fig:exea-MAL-0.5}}\\
\caption{Example \ref{ex1}: The effect using  $\tilde L_\text{bdry}({\bs{\theta}})$} 
\label{Fig:4_revision}
\end{figure} 

\begin{remark}[How does $V_{\bs{\eta}}$ affect the  method's performance?]
In the proof, we require $V_{\bs{\eta}}$ to be rich enough to satisfy the stability condition. In this example, we tested multiple configurations for the $V_{\bs{\eta}}$ network, experimenting with different hidden layers and varying numbers of neurons. The results are all consistent with  \Cref{Fig:1}. This indicates that the model is robust with $V_{\bs{\eta}}$ for this example.
\end{remark}
%



\subsection{Stationary PDE problems}\label{subsec:statPDE}

\begin{ex}\label{ex2}
We now test the following high-dimensional problem as in \cite{ZBYZ20}.
\[
	\begin{cases}
		- \triangle  u(\bs{x}) = f & \bs{x} \in \Omega\\
		u(\bs{x}) = g (\bs{x})&  \bs{x} \mbox{ on } \partial \Omega,
	\end{cases}
\]
where the true solution renders
$u(\bs{x}) = \sum\limits_{i=1}^{d} \sin\left(\dfrac{\pi}{2}x_{i}\right)$.

\end{ex}

Observe that the boundary condition on the plane $x_{1}=0$ and $x_{1}=1$ now serves as the initial and terminal conditions in the pseudo-time XNODE model. Subsequently, we adjust the initial loss and introduce the terminal loss as:
\begin{equation}
\begin{split}
&L_\text{init}(\theta)= ||u_\theta(x_{1} = 0 )- g (x_{1} = 0)||_{L^2(\Omega(0))},\\
&L_\text{last}(\theta):= ||u_\theta(x_{1} = 1 )- g (x_{1} = 1)||_{L^2(\Omega(1))},
\end{split}
\end{equation}
where $\Omega(t) := \{ \bs{x} \in \Omega, \; x_{1} = t\}$.
We shall utilize the following loss functions for the pseudo-time XNODE model.
\begin{equation}\label{loss-b}
\tilde L_{\text{wan},\text{cwan},\text{scwan}}(\theta, \eta) =L_{\text{wan},\text{cwan},\text{scwan}}(\theta, \eta) + \gamma L_\text{last}(\theta).
\end{equation}

We experimented with testing the pseudo-time XNODE model with $d=5$, using the same parameters as in \cref{tab:1} except for the penalty parameters. A grid search was performed to tune the hyperparameters $\alpha$, $\beta$, and $\gamma$, which were restricted to the range $[10, 10^9]$. A optimal values found were $\alpha = \gamma =10^{5}$ and $\beta = 10^{7}$.

We established a stopping criteria that ensured the relative training error was below $1\%$ or the maximum iteration number less than 300. Each training step takes about $8$ seconds.

We have analyzed the loss functions $\tilde L_{\text{wan}}$, $\tilde L_{\text{scwan}}$ with $\gamma_{d}=0.5$, and $\tilde L_{\text{cwan}}$ with $\gamma_{d}=0.001$, as described in \cref{loss-b}. We conducted tests on both the recursive and non-recursive models for each setting, and the outcomes are displayed in \Cref{Fig:ex2}. Each sub-figure in \Cref{Fig:ex2} showcases the results of three consecutive experiments that were initialized randomly. The recursive and non-recursive models are utilized in the left and right columns respectively.

For $\tilde L_{\text{wan}}$, the recursive model in \Cref{fig:ex2-wan-r} reached the stopping criteria at steps $169, 98$, and $91$ (with an average of $120$). Meanwhile, the non-recursive model in \Cref{fig:ex2-wan} reached the stopping criteria at steps $277, 119$, and $93$ (with an average of $163$), based on three experiments for each.

For $\tilde L_{\text{cwan}}$, the recursive model in \Cref{fig:ex2-cwan-r} reached the stopping criteria at steps $237, 174$, and $200$ (with an average of $203$). Meanwhile, the non-recursive model in \Cref{fig:ex2-cwan} reached the stopping criteria at steps $293, 149$, and $153$ (with an average of $198$), based on three experiments for each.

For $\tilde L_{\text{scwan}}$, the recursive model in \Cref{fig:ex2-scwan-r} reached the stopping criteria at steps $172, 144$, and $117$ (with an average of $144$). Meanwhile, the non-recursive model in \Cref{fig:ex2-scwan} reached the stopping criteria at steps $151, 132$, and $115$ (with an average of $132$), based on three experiments for each.

\begin{figure}[!ht]
    \centering
            \subfloat[XNODE(R) + $\tilde L_{\text{wan}}$]{\includegraphics[width=0.48\textwidth]{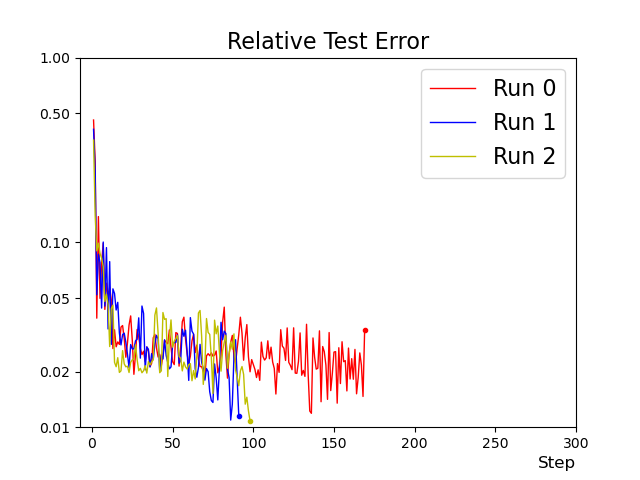}
         \label{fig:ex2-wan-r}}
                 \subfloat[XNODE + $\tilde L_{\text{wan}}$]{\includegraphics[width=0.48\textwidth]{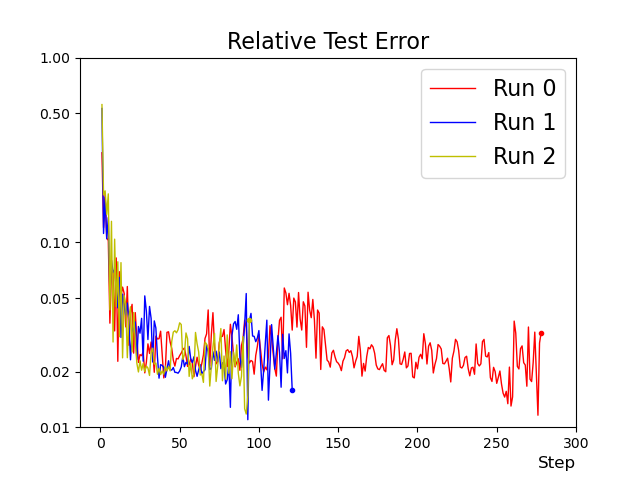}
         \label{fig:ex2-wan}}
         \\
                 \subfloat[XNODE(R) + $\tilde L_{\text{cwan}} (\gamma_{d}=10^{-3})$]{\includegraphics[width=0.48\textwidth]{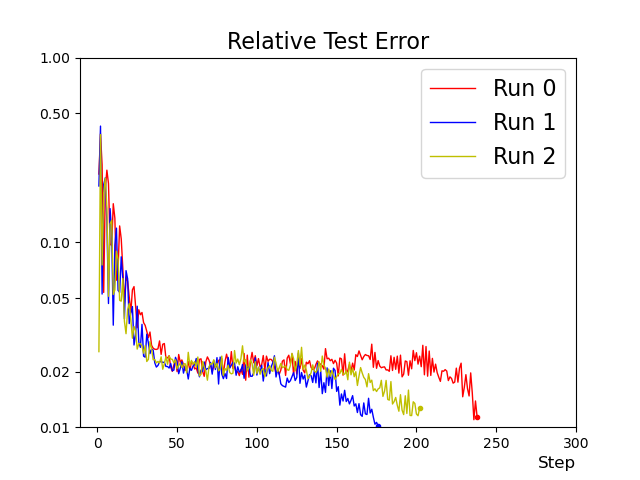}
         \label{fig:ex2-cwan-r}}
         \subfloat[XNODE+ $ \tilde L_{\text{cwan}} (\gamma_{d}=10^{-3})$]{\includegraphics[width=0.48\textwidth]
         {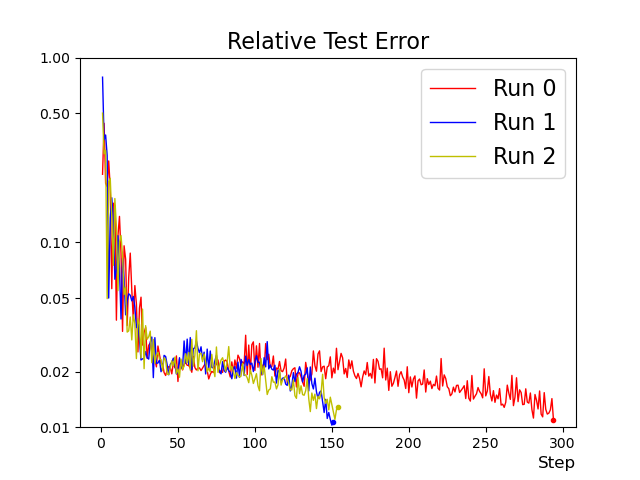}
         \label{fig:ex2-cwan}}\\
        \subfloat[XNODE+ $\tilde L_{\text{scwan}} (\gamma_{d}=0.5)$]{\includegraphics[width=0.48\textwidth]
        {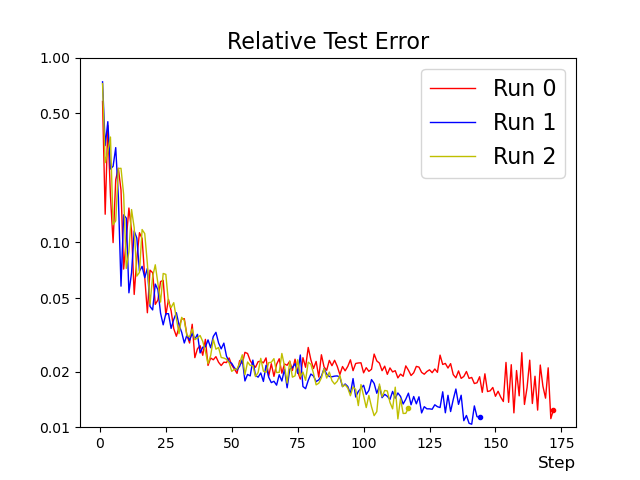}
         \label{fig:ex2-scwan-r}}
         \subfloat[XNODE(R)+ $ \tilde L_{\text{Scwan}} (\gamma_{d}=0.5)$]{\includegraphics[width=0.48\textwidth]
         {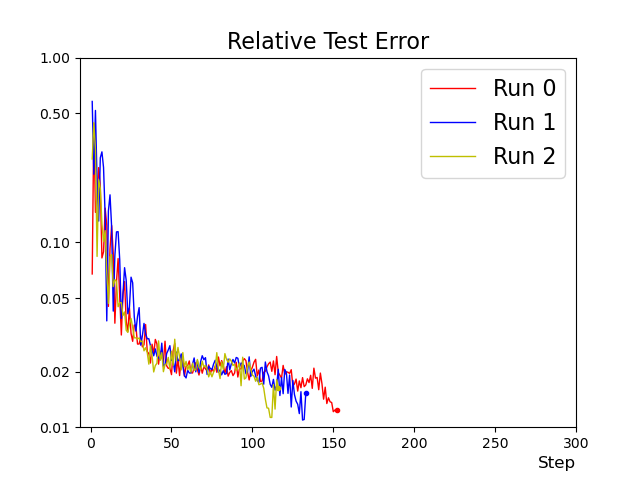}
         \label{fig:ex2-scwan}}
\caption{Example \ref{ex2}. Relative $L^{2}$ Error versus Step using pseudo-time XNODE} \label{Fig:ex2}
\end{figure}

In all XNODE experiments, the relative error quickly reached the $2\%$ threshold within the first $35$ iterations.
Although the relative error oscillations generated by $\tilde L_{\text{wan}}$ are still greater than those of $\tilde L_{\text{cwan}}$ and $\tilde L_{\text{scwan}}$, it is worth noting that, in this particular case, the stopping criteria was achieved with slightly fewer iterations on average. {We believe this faster convergence takes place thanks to the Poisson type PDE used in this example. For the Poisson problem, it is easy to see that $\|\cdot\|_{op}$ is the most natural norm to minimize.}
It has also been observed that the recursive model's performance is almost comparable to that of the non-recursive models in this example.

\begin{ex}\label{ex3}
\begin{equation}
-\nabla \cdot (a(x) \nabla u) + \dfrac{1}{2} |\nabla u|^{2} = f(x), \mbox{ in } \Omega = [0,1]^{d},\quad
u(x) = g(x) \mbox{on } \partial \Omega.
\end{equation}
where $a(x) = 1 + \|x\|^{2}$. The true solution $u(x) = \sin(0.5\pi x_{1}^{2}+ 0.5 x_{2}^{2})$.
\end{ex}
In this problem, the non-linear term $\dfrac{1}{2} |\nabla u|^{2}$ presents a significant challenge. We will use the hyper-parameter set from \Cref{{tab:2}} in all numerical tests with the pseudo-time XNODE model. Our objective in this example is to evaluate the performance of the pseudo-time XNODE model using various loss functions. We established the stopping criteria for the maximum number of iterations to be less than $600$.
 The duration of each iteration is approximately $8.5$ seconds.

We first test the loss function $\tilde L_{\text{wan}}$ by conducting three consecutive experiments with random initialization. The results are presented in \Cref{Fig:ex3-wan}. The left/ right figure in \Cref{Fig:ex3-wan} shows the relationship between the number of steps and the $L^{2}$ relative error/ minimal $L^{2}$ relative error based on test sets. After the stopping criteria have been met, the minimal relative training $L^{2}$ error is $0.024$, $0.056$, and $0.023$, respectively. We have observed a slower convergence rate compared to previous examples, which can be attributed to the more challenging non-linear term.

\begin{table}[!ht]\label{tab:2}
\small{
\begin{tabular}{|l|l|l|l|l|l|l|l|l|}
\hline
 $d$ &$N_r$  & $N_b$ &$n_{T}$  & $K_{u}$  & $K_{\phi}$  & $\alpha$  & $\beta$ & $\gamma$  \\ \hline
 $5$  &  $4000$ & $4000$    & $20$          &  $2$            &  $1$ & $6\times10^{7}$ & $12\times10^{7}$& $12\times10^{7}$   \\ \hline
 $\epsilon$ & $l_{\theta}$  & $l_{\eta}$  & $u_{\text{layers}}$  & $u_{\text{hid-dim}1}$  &$u_{\text{hid-dim}2}$  	 
  &$v_{\text{layers}}$  
  & $v_{\text{hid-dim}}$ &  \\ \hline
 $10^{-2}$&  $.015$& $.03$  & $12$  & $20$  & $10$  & $9$  &  $50$ & \\ \hline
\end{tabular}
\caption{Hyper parameter setting for Example \ref{ex3}}
}
\end{table}
\begin{figure}[!h]
    \centering
      \includegraphics[width=0.48\textwidth]{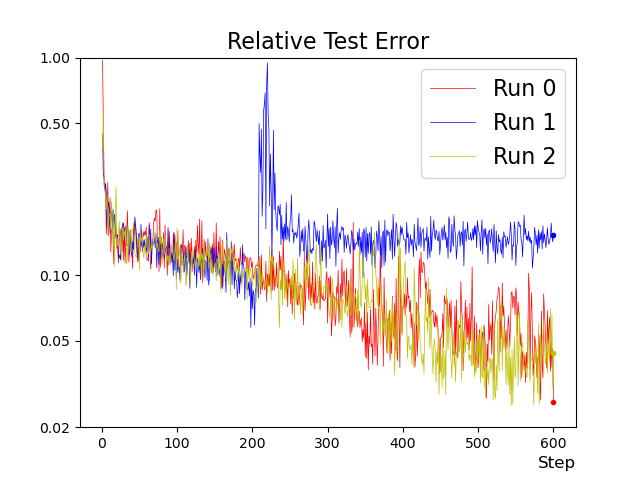} 
      \includegraphics[width=0.48\textwidth]{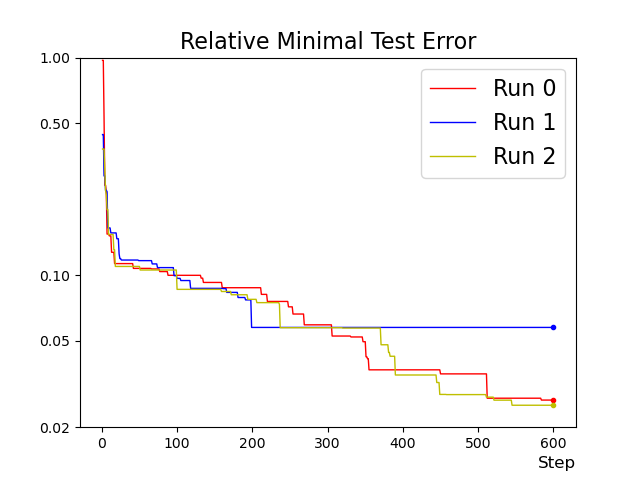}
\caption{Example \ref{ex3}. Pseudo time XNODE + $\tilde L_{\text{wan}}$} \label{Fig:ex3-wan}
\end{figure} 

The results for $\tilde L_{\text{cwan}}$ and $\tilde L_{\text{scwan}}$ both with $\gamma_{d}=0.5$ are provided in \Cref{Fig:ex3-AL} and \Cref{Fig:ex3-MAL}, respectively.
For all three experiments, the minimal relative error calculated from $\tilde L_{\text{cwan}}$ was $0.013, 0.016$, and $0.016$. Meanwhile, the minimal relative error calculated from $\tilde L_{\text{scwan}}$ for the same experiments were $0.011, 0.012$, and $0.016$.
Therefore, in comparison to \Cref{Fig:ex3-wan}, using $\tilde L_{\text{cwan}}$ and $\tilde L_{\text{scwan}}$ provides a slight advantage over 
$\tilde L_{\text{wan}}$.
Due to the highly nonlinear nature of this problem, we believe that a more refined approach to tuning the hyperparameters is necessary to achieve greater accuracy in the results. We will consider this as a future work.

\begin{figure}[!h]
    \centering
      \includegraphics[width=0.48\textwidth]{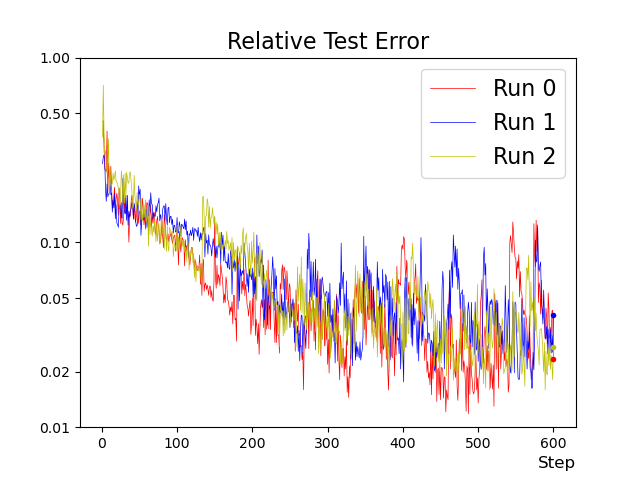} 
      \includegraphics[width=0.48\textwidth]{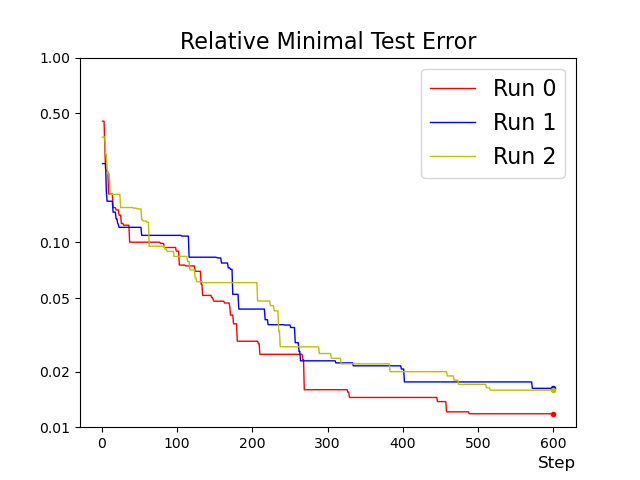}
\caption{Example \ref{ex3}. Pseudo time XNODE + $\tilde L_{\text{cwan}}$ ($\gamma_{d} =0.5$)} \label{Fig:ex3-AL}
\end{figure} 

\begin{figure}[!h]
    \centering
      \includegraphics[width=0.48\textwidth]{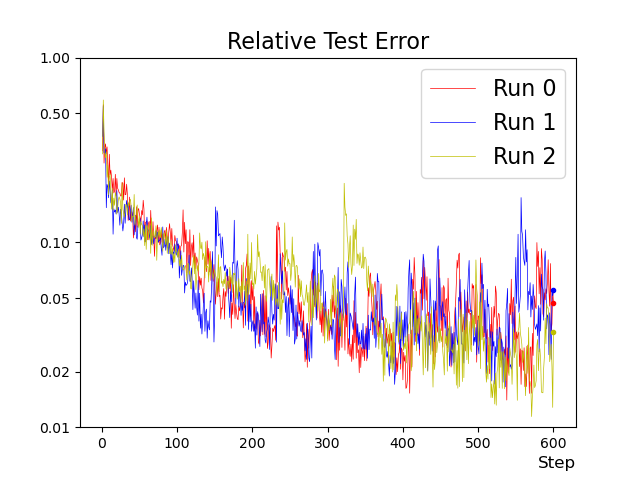} 
      \includegraphics[width=0.48\textwidth]{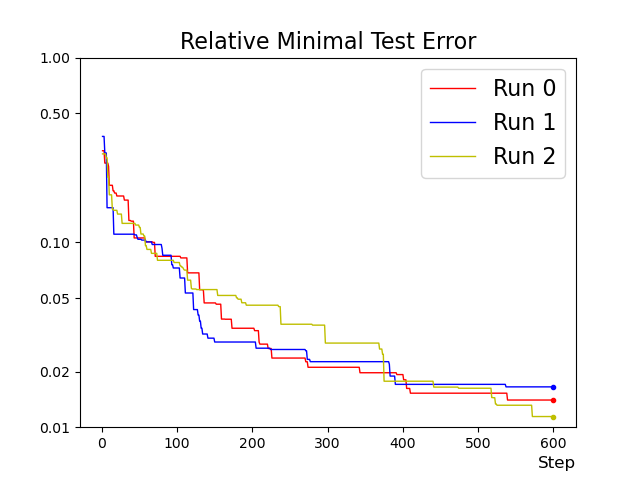}
\caption{Example \ref{ex3}. Pseudo time XNODE + $\tilde L_{\text{scwan}}$ ($\gamma_{d} =0.5$)} \label{Fig:ex3-MAL}
\end{figure}

\appendix
\section{Model set up XNODE-WAN Algorithm}\label{appen:A}
 The hyper-parameters for the neural networks are explained in the following table.
\begin{table}[h]
\footnotesize{
\centering
\resizebox{\columnwidth}{!}{%
\begin{tabular}{l|l}
\hline
         Notation  & Meaning \\
         \hline
$d$ &  dimension for the physical domain (not including the time domain)\\
$N_r$      & Number of sampled collocation points of the spatial domain \\
$N_b$ &  Number of sampled collocation points of the spatial domain boundary  \\
$n_T$ &  Number of sampled time partition \\
$K_u$ &  Inner iteration to update weak solution $u_\theta$ or $u_\theta$ \\
$K_\phi$ &  Inner iteration to update test function $\phi_\eta$  \\
$\alpha$ &   Weight parameter of boundary loss ${L}_\text{bdry}$\\
$\gamma$ &  Weight parameter of initial and terminal losses ${L}_\text{init}$ and ${L}_\text{last}$\\
$\epsilon$ & relative error tolerance\\
$l_\theta$  & Learning rate for the primal network \\
$l_\eta$ & Learning rate for network parameter $\eta$ of test function $\phi_\eta$\\
$u_{\text{layers}}$ & The number of hidden layers for $\mathcal{N}^{\text{vec}}_{\theta_2}$\\
$u_{\text{hid-dim}1}$ & intermediate and output dimension for $\mathcal{N}^{\text{init}}_{\theta_{1}}$, input dimension for $\mathcal{N}^{\text{vec}}_{\theta_2}$\\
$u_{\text{hid-dim}2}$ &  intermediate dimension for $\mathcal{N}^{\text{vec}}_{\theta_2}$\\
$v_{\text{layers}}$ & The number of hidden layers for $V_{\bs{\eta}}$\\
$v_{\text{hidd-dim}}$ & intermediate and output dimension for $V_{\bs{\eta}}$.\\
\hline
\end{tabular}}
\caption{List of hyperparameters.}
\label{table2}
}
\end{table}
For $V_{\bs{\eta}}$, we use a classical DNN network. The activation is set to be Tanh for the last hidden layer and ReLU for other hidden layers. Note that there is no activation function for the output layer.
$\mathcal{N}^{\text{init}}_{\bs{\theta}_{1}}$ has one input layer, one hidden layer, and one output layer. The activation function after both the input and hidden layer is ReLU.
$\mathcal{N}^{\text{vec}}_{\bs{\theta}_{2}}$ has $u_{\text{layers}}$ of hidden dimensions and $\text{Tanh}$ as the activation function.

\bibliographystyle{siamplain}
\bibliography{references}

\end{document}


\maketitle

%
%
%
%
%
Below we provide the pseudo algorithms for the XNODE and pseudo-time XNODE methods in approximating (5.3).
\begin{algorithm}[H]
	{
		\caption{XNODE Network for approximating (5.3)}
	\label{alg:XNODE}
	\textbf{Input:}  A sample point $\bs{x} \in \bar{\Omega}$ and $\Pi_{T} = (t_{i})_{i =0}^{n_T}$ time partition of $[0, T]$.\\
	\textbf{Algorithm:}
	\begin{algorithmic}[1]
		\STATE compute $o_{\bs{x}}(\Pi_{T})$  by approximating (5.3) using e.g., forward Euler method.
		\RETURN  $o_{\bs{x}}(\Pi_{T})$   \# \texttt{ $o_{\bs{x}}(\Pi_{T})_{1 \times n_{T}}$}
	\end{algorithmic}
}
\end{algorithm}

{\begin{algorithm}[H] 
		\small{
			\caption{Pseudo-time XNODE Network for approximating (5.3)} 
		\label{alg:sudoXNODE}
		\textbf{Input:}  A sample point $\tilde{\bs{x}} \in [0,1]^{d-1}$ and $\Pi_{T} = \{x_{1,i}\}_{i =0}^{n_T}$ pseudo time partition of $[0, 1]$.\\
		\textbf{Algorithm:}
		\begin{algorithmic}[1] 
			\STATE compute $o_{\tilde{\bs{x}}}(\Pi_{T})$  by approximating (5.3) using e.g., the forward Euler method.
			\RETURN output  $o_{\tilde{\bs{x}}}(\Pi_{T})$  \# \texttt{ $o_{\tilde{\bs{x}}}(\Pi_{T})_{1 \times n_{T}}$}
	\end{algorithmic}}
\end{algorithm}}
{We now provide the complete pseudo algorithm for the WAN method with XNODE and pseudo-time XNODE networks in \Cref{XNODE_alg}.
	\begin{algorithm}[ht!]
			\caption{XNODE/ Pseudo-time XNODE Algorithm}\label{XNODE_alg}
			\textbf{Input:} domain $\mathcal{D}$, tolerance $\epsilon>0$, $N_r/N_b$: number of sampled points on the domain/boundary condition,  $K_u/K_\phi$: number of solution/adversarial network parameter updates per iteration, $\alpha/\gamma \in \mathbb{R}^+:$ the weight of boundary/initial loss, $N_{max}$: the maximal iteration number,  and the hyperparameters for all DNNs involved
			\begin{algorithmic}[1]
				\STATE{\textbf{Initialise:} }
				\STATE{ \hspace{1cm}The network parameters $\theta$ and $\eta$ for $u_\theta,\phi_\eta:\mathcal{D} \to \mathbb{R}$}, respectively\\ 
			\STATE{\hspace{1cm}Generate points sets $S^{r}$ and $S^{b}$} 
			\texttt{\#random point sampling}
			\STATE{\hspace{1cm}Solve for $u_{\bs{\theta}}$ using  \Cref{alg:XNODE} or \Cref{alg:sudoXNODE}}
			\STATE{$Ite = 0$} \texttt{\# the number of iteration}
			\WHILE{$\dfrac{\|u - u_{\bs{\theta}}\|_{\mathcal{D}}}{\|u\|_{\mathcal{D}}} \ge \epsilon$ and $Ite < N_{max}$}  
		\STATE{\# \texttt{update weak solution network parameters}}
		\FOR{$k = 1,...,K_u$}
		\STATE{compute $\nabla_\theta  L(\theta, \eta)$;} 
		\STATE{update $\theta\gets\theta- l_\theta \nabla_\theta  L(\theta, \eta)$;} 
		\ENDFOR
		\STATE{\# \texttt{update test network parameters}}
		\FOR{$k = 1,...,K_\phi$}
		\STATE{compute $\nabla_\eta  L(\theta, \eta)$;}
		\STATE{update \small{$\eta\gets\eta+l_\eta \nabla_\eta  L(\theta, \eta)$;}
		}
		\ENDFOR
		\STATE{$Ite = Ite +1$}
		\STATE{Regenerate points sets $S^{r}$ and $S^{b}$;}\\  
		\STATE{Solve for $u_{\bs{\theta}}$ using  \Cref{alg:XNODE} or \Cref{alg:sudoXNODE}}
		
		\ENDWHILE\\
		\textbf{Output:} the weak solution $u_\theta:\mathcal{D} \to\mathbb{R}$
	\end{algorithmic}
\end{algorithm}}

\section{Extension to nonlinear problems}

We briefly sketch how, under suitable assumptions, the theoretical results presented in Section 2.1 can be extended to a nonlinear problem. Let $A: W \to W'$ denote, throughout this section, a non linear operator satisfying the following assumptions:
\begin{enumerate}
	\item $A$ is monotone: there exists $\alpha > 0$ such that for all $w, w' \in W$ it holds that
	\[
	\langle A(w) - A(w'), w - w' \rangle \geq \alpha \| w - w' \|_W^2;
	\]
	\item $A$ is Lipschitz continuous: there exists $M$ such that for all $w, w' \in W$ it holds that
	\[
	\| A(w) - A(w') \|_{W'} \leq M \| w - w' \|_{W};
	\]
	\item $0 \in \mathrm{Range}(A)$: there exists $u \in W$ such that $A(u) = 0$. \label{existence}
\end{enumerate}

We consider the nonlinear equation 
\begin{equation}\label{nonlin}
A(u) = 0,
\end{equation}
which, given $W_{\bs\theta}$ and $V_{\bs\eta}$ satisfying the assumptions of Section 2.1, we approximate with $u_{\bs\theta}^*$ defined as
\[
u^*_{\bs\theta} = \mathop{\mathrm{argmin}}_{w_{\bs\theta} \in W_{\bs\theta}} \, \sup_{v_{\bs\eta} \in V_{\bs\eta}} \frac{\langle A(w_{\bs\theta}),v_{\bs\eta} \rangle}{\| v_{\bs\eta} \|_{W}}.
\]
The analysis of the discrete problem can be mostly carried out as in the linear case. 
We observe that the supremum on the right hand side is necessarily non negative, and therefore it admits an infimum $\sigma^*$ for the minimization in $w_{\bs\theta}$. Let $\{ w_{\bs\theta}^n\}$ denote a minimizing sequence, that is, 
$\sup_{v_{\bs\eta} \in V_{\bs\eta}} \frac{\langle A(w_{\bs\theta}^n),v_{\bs\eta} \rangle}{\| v_{\bs\eta} \|_{W}}$
converges to $\sigma^*$ as $n$ goes to $\infty$. Letting $\widehat w_{\bs\theta}$ denote an arbitrary but fixed element of $W_{\bs\theta}$ we can write
\begin{equation*}
	\begin{split}
\| w^n_{\bs\theta} \|^2_{W} &\leq 2 \| w^n_{\bs\theta} -  \widehat w_{\bs\theta}\|_{W}^2 + 2 \| \widehat w_{\bs\theta} \|_{W} ^2 \\&\leq \frac 2 \alpha \langle A(w_{\bs\theta}^n) - A(\widehat w_{\bs\theta}), w_{\bs\theta}^n - \widehat w_{\bs\theta} \rangle + 2 \| \widehat w_{\bs\theta} \|_{W} ^2 \\
&\leq \frac {2M} \alpha (\|  A(w_{\bs\theta}^n) \|_{W'} + \| A(\widehat w_{\bs\theta}) \|_{W'} ) (\| w^n_{\bs\theta} \|_W + \| \widehat w_{\bs\theta} \|_W) + 2 \| \widehat w_{\bs\theta} \|_W^2\\
& \leq \frac {4 M^2} {\alpha^2} \left(\sup_{v_{\bs\eta} \in V_{\bs\eta}} \frac {\langle A(w^n_{\bs\theta}), v_{\bs\eta} \rangle } {\| v_{\bs\eta} \|_{W}}  + \| A (\widehat w_{\bs\theta}) \|_{W'}\right)^2 + \frac 1 2 \| w^n_{\bs\theta} \|^2_{W} + \frac 5 2 \| \widehat w_{\bs\theta} \|_W^2.
\end{split}
\end{equation*}
Subtracting $\| w^n_{\bs\theta} \|_{W} / 2$ from both sides we see that
\begin{equation*}
	\frac 1 2 \| w^n_{\bs\theta} \|^2_{W}  \leq \frac {4 M^2} {\alpha^2} \left(\sup_{v_{\bs\eta} \in V_{\bs\eta}} \frac {\langle A(w^n_{\bs\theta}), v_{\bs\eta} \rangle } {\| v_{\bs\eta} \|_{W}}  + \| A (\widehat w_{\bs\theta}) \|_{W'}\right)^2 + \frac 5 2 \| \widehat w_{\bs\theta} \|_W^2,
\end{equation*}
from which we see that the sequence $\{ w_{\bs\theta}^n\}$ is bounded in $W$. Then,  as in Section 2.1, it admits a weakly convergent subsequence $\{\tun\}$, with weak limit $\utstar \in \mathrm{cl}^{seq}_w(W_{\bs\theta})$.
We claim that $\utstar$ satisfy the approximation bound given in Lemma 2. Indeed, letting $w_{\bs\theta}$ be an arbitrary element of $W_{\bs\theta}$ we have
\begin{equation*}
		\| \tun - w_{\bs\theta} \|_{W}  \leq \frac 1 \alpha \frac {\langle A(\tun) - A(w_{\bs\theta}), \tun - w_{\bs \theta} \rangle }{\| \tun - w_{\bs\theta} \|_{W} } \leq \frac 1 \alpha \sup_{v_{\bs\eta}\in V_{\bs\eta}}  \frac {\langle A(\tun) - A(w_{\bs\theta}), v_{\bs\eta} \rangle }{\| v_{\bs\eta} \|_{W} }.
\end{equation*}

Then, as the norm of the weak limit of a weakly convergence sequence is bounded from above by the limit of the norms of the sequence elements, we can write
\begin{equation*}
	\begin{split}
		\| w_{\bs\theta} - \utstar \|_W & \leq \lim_{n \to \infty} \| w_{\bs\theta} - \tun \|_{W} 
		 \leq \frac 1 \alpha \lim_{n \to \infty} \sup_{v_{\bs\eta} \in V_{\bs\eta}} \frac{
	\langle A(\tun) - A(w_{\bs\theta}) , v_{\bs\eta}	\rangle }{\| v_{\bs\eta} \|_{W}} 
		\\& \leq \frac 1 \alpha \left(  \lim_{n \to \infty}
	\sup_{v_{\bs\eta} \in V_{\bs\eta}} \frac{
		\langle A( \tun )  , v_{\bs\eta}	\rangle }{\| v_{\bs\eta} \|_{W}} + \sup_{v_{\bs\eta} \in V_{\bs\eta}}   \frac{
		\langle  A(w_{\bs\theta}) , v_{\bs\eta}	\rangle }{\| v_{\bs\eta} \|_{W}}
	\right) 
	\\& \leq \frac 1 \alpha \left(  \inf_{w'_{\bs\theta} \in W_{\bs\theta}}
	\sup_{v_{\bs\eta} \in V_{\bs\eta}} \frac{
		\langle A( w'_{\bs\theta} )  , v_{\bs\eta}	\rangle }{\| v_{\bs\eta} \|_{W}} + \sup_{v_{\bs\eta} \in V_{\bs\eta}}   \frac{
			\langle  A(w_{\bs\theta}) , v_{\bs\eta}	\rangle }{\| v_{\bs\eta} \|_{W}}
	\right) \\
&\leq \frac 2 \alpha  \sup_{v_{\bs\eta} \in V_{\bs\eta}}   \frac{
		\langle  A(w_{\bs\theta}) , v_{\bs\eta}	\rangle }{\| v_{\bs\eta} \|_{W}} = \frac 2 \alpha  \sup_{v_{\bs\eta} \in V_{\bs\eta}}   \frac{
		\langle  A(w_{\bs\theta})  - A(u), v_{\bs\eta}	\rangle }{\| v_{\bs\eta} \|_{W}}\\ & \leq 2 \frac M \alpha \| w_{\bs\theta} - u \|_{W}.
		\end{split}
\end{equation*}

Using a triangle inequality and exploiting the arbitrariness of $w_{\bs\theta}$ we obtain that
\[
\| u - \utstar \|_{W} \leq \left( 1 + 2 \frac M \alpha\right) \inf_{w_{\bs\theta}} \| u - w_{\bs\theta} \|_W.
\]

%% file: ex_shared.tex
\usepackage{lipsum}
\usepackage{amsfonts}
\usepackage{subfig}
\usepackage{graphicx}
\usepackage{epstopdf}
\usepackage{algorithmic}
\usepackage{amsmath,amssymb}
\usepackage{mathabx}
\usepackage{xfrac}
\ifpdf
  \DeclareGraphicsExtensions{.eps,.pdf,.png,.jpg}
\else
  \DeclareGraphicsExtensions{.eps}
\fi

\newsiamremark{remark}{Remark}
\newsiamremark{hypothesis}{Hypothesis}
\crefname{hypothesis}{Hypothesis}{Hypotheses}
\newsiamthm{claim}{Claim}

\headers{WAN discretization of PDEs
}{S. Bertoluzza, E. Burman, and C. He}

\title{WAN discretization of PDEs: best approximation, stabilization and essential boundary conditions\thanks{Submitted to SIAM Journal of Scientific Computing.
\funding{EB was partially funded by EPSRC grants EP/W007460/1, EP/V050400/1 and EP/T033126/1. SB was partially funded by PRIN grant 20204LN5N5 AdPoly}}}

\title{WAN discretization of PDEs: best approximation, stabilization and essential boundary conditions\thanks{Submitted to SIAM Journal of Scientific Computing.
\funding{EB was partially funded by EPSRC grants EP/W007460/1, EP/V050400/1 and EP/T033126/1.}}}
\author{Silvia Bertoluzza\thanks{Istituto di Matematica Applicata e Tecnologie Informatiche, CNR, Italy,
  (\email{silvia.bertoluzza@imati.cnr.it}).}
\and Erik Burman\thanks{Department of  Mathematics, University College London, UK, 
  (\email{e.burman@ucl.ac.uk}).}
\and Cuiyu He \thanks{Department of Mathematics, Oklahoma State University, 401 Stillwater, OK, 74078 (\email{cuiyu.he@okstate.edu}).// Department of Mathematics, University of Georgia, Athens, GA 30602 (\email{cuiyu.he@uga.edu})} }

\usepackage{amsopn}


%% file: main.bbl
\begin{thebibliography}{10}

\bibitem{BYZZ20}
{\sc G.~Bao, X.~Ye, Y.~Zang, and H.~Zhou}, {\em Numerical solution of inverse
  problems by weak adversarial networks}, Inverse Problems, 36 (2020),
  pp.~115003, 31, \url{https://doi.org/10.1088/1361-6420/abb447},
  \url{https://doi.org/10.1088/1361-6420/abb447}.

\bibitem{BFMO21}
{\sc E.~Burman, A.~Feizmohammadi, A.~M\"{u}nch, and L.~Oksanen}, {\em Space
  time stabilized finite element methods for a unique continuation problem
  subject to the wave equation}, ESAIM Math. Model. Numer. Anal., 55 (2021),
  pp.~S969--S991, \url{https://doi.org/10.1051/m2an/2020062},
  \url{https://doi.org/10.1051/m2an/2020062}.

\bibitem{BO18}
{\sc E.~Burman and L.~Oksanen}, {\em Data assimilation for the heat equation
  using stabilized finite element methods}, Numer. Math., 139 (2018),
  pp.~505--528, \url{https://doi.org/10.1007/s00211-018-0949-3},
  \url{https://doi.org/10.1007/s00211-018-0949-3}.

\bibitem{deshmukh2016neural}
{\sc A.~N. Deshmukh, M.~Deo, P.~K. Bhaskaran, T.~B. Nair, and K.~Sandhya}, {\em
  Neural-network-based data assimilation to improve numerical ocean wave
  forecast}, IEEE Journal of Oceanic Engineering, 41 (2016), pp.~944--953.

\bibitem{DP94}
{\sc M.~W. M.~G. Dissanayake and N.~Phan-Thien}, {\em Neural-network-based
  approximations for solving partial differential equations}, Communications in
  Numerical Methods in Engineering, 10 (1994), pp.~195--201,
  \url{https://doi.org/https://doi.org/10.1002/cnm.1640100303},
  \url{https://onlinelibrary.wiley.com/doi/abs/10.1002/cnm.1640100303},
  \url{https://arxiv.org/abs/https://onlinelibrary.wiley.com/doi/pdf/10.1002/cnm.1640100303}.

\bibitem{duan2021analysis}
{\sc C.~Duan, Y.~Jiao, Y.~Lai, X.~Lu, Q.~Quan, and J.~Z. Yang}, {\em Analysis
  of deep ritz methods for laplace equations with dirichlet boundary
  conditions}, 2021, \url{https://arxiv.org/abs/2111.02009}.

\bibitem{DL20}
{\sc M.~Duprez and A.~Lozinski}, {\em {$\phi$}-{FEM}: a finite element method
  on domains defined by level-sets}, SIAM J. Numer. Anal., 58 (2020),
  pp.~1008--1028, \url{https://doi.org/10.1137/19M1248947},
  \url{https://doi.org/10.1137/19M1248947}.

\bibitem{ern2004theory}
{\sc A.~Ern and J.-L. Guermond}, {\em Theory and practice of finite elements},
  vol.~159, Springer, 2004.

\bibitem{guhring2020error}
{\sc I.~G{\"u}hring, G.~Kutyniok, and P.~Petersen}, {\em Error bounds for
  approximations with deep relu neural networks in w s, p norms}, Analysis and
  Applications, 18 (2020), pp.~803--859.

\bibitem{guhring2020expressivity}
{\sc I.~G{\"u}hring, M.~Raslan, and G.~Kutyniok}, {\em Expressivity of deep
  neural networks}, Cambridge University Press, 2022.

\bibitem{he2016deep}
{\sc K.~He, X.~Zhang, S.~Ren, and J.~Sun}, {\em Deep residual learning for
  image recognition}, in Proceedings of the IEEE conference on computer vision
  and pattern recognition, 2016, pp.~770--778.

\bibitem{hong2021priori}
{\sc Q.~Hong, J.~W. Siegel, and J.~Xu}, {\em A priori analysis of stable neural
  network solutions to numerical pdes}, arXiv preprint arXiv:2104.02903,
  (2021).

\bibitem{karumuri2020simulator}
{\sc S.~Karumuri, R.~Tripathy, I.~Bilionis, and J.~Panchal}, {\em
  Simulator-free solution of high-dimensional stochastic elliptic partial
  differential equations using deep neural networks}, Journal of Computational
  Physics, 404 (2020), p.~109120.

\bibitem{mahan2021nonclosedness}
{\sc S.~Mahan, E.~J. King, and A.~Cloninger}, {\em Nonclosedness of sets of
  neural networks in sobolev spaces}, Neural Networks, 137 (2021), pp.~85--96.

\bibitem{moeini2012wave}
{\sc M.~H. Moeini, A.~Etemad-Shahidi, V.~Chegini, and I.~Rahmani}, {\em Wave
  data assimilation using a hybrid approach in the persian gulf}, Ocean
  Dynamics, 62 (2012), pp.~785--797.

\bibitem{oliva2021fast}
{\sc P.~V. Oliva, Y.~Wu, C.~He, and H.~Ni}, {\em Towards fast weak adversarial
  training to solve high dimensional parabolic partial differential equations
  using xnode-wan}, 2021, \url{https://arxiv.org/abs/2110.07812}.

\bibitem{ostrovskii2001weak}
{\sc M.~Ostrovskii}, {\em Weak* sequential closures in {B}anach space theory
  and their applications}, 2001,
  \url{https://arxiv.org/abs/arXiv:math/0203139}.

\bibitem{petersen2021topological}
{\sc P.~Petersen, M.~Raslan, and F.~Voigtlaender}, {\em Topological properties
  of the set of functions generated by neural networks of fixed size},
  Foundations of computational mathematics, 21 (2021), pp.~375--444.

\bibitem{raissi2019physics}
{\sc M.~Raissi, P.~Perdikaris, and G.~E. Karniadakis}, {\em Physics-informed
  neural networks: A deep learning framework for solving forward and inverse
  problems involving nonlinear partial differential equations}, Journal of
  Computational Physics, 378 (2019), pp.~686--707.

\bibitem{stein70}
{\sc E.~M. Stein}, {\em Singular integrals and differentiability properties of
  functions}, Princeton Mathematical Series, No. 30, Princeton University
  Press, Princeton, N.J., 1970.

\bibitem{sukumar2022exact}
{\sc N.~Sukumar and A.~Srivastava}, {\em Exact imposition of boundary
  conditions with distance functions in physics-informed deep neural networks},
  Computer Methods in Applied Mechanics and Engineering, 389 (2022), p.~114333.

\bibitem{ZBYZ20}
{\sc Y.~Zang, G.~Bao, X.~Ye, and H.~Zhou}, {\em Weak adversarial networks for
  high-dimensional partial differential equations}, J. Comput. Phys., 411
  (2020), pp.~109409, 14, \url{https://doi.org/10.1016/j.jcp.2020.109409},
  \url{https://doi.org/10.1016/j.jcp.2020.109409}.

\bibitem{zeng2022adaptive}
{\sc S.~Zeng, Z.~Zhang, and Q.~Zou}, {\em Adaptive deep neural networks methods
  for high-dimensional partial differential equations}, Journal of
  Computational Physics, 463 (2022), p.~111232.

\end{thebibliography}
